\documentclass[ot]{birkmult}
\usepackage{amssymb}
\usepackage{amsmath}
\usepackage{graphicx,color}
\usepackage{latexsym}
\usepackage{amsfonts}

\newtheorem{theorem}{Theorem}[section]
\newtheorem{lemma}[theorem]{Lemma}

\theoremstyle{definition}
\newtheorem{definition}[theorem]{Definition}

\newtheorem{remark}[theorem]{Remark}

\numberwithin{equation}{section}

\newcommand{\calA}{{\mathcal A}}

\newcommand{\calH}{{\mathcal H}}

\begin{document}
\title[]{The Schur Transformation for Nevanlinna \\ Functions:
Operator Representations, \\ Resolvent Matrices, and Orthogonal \\
Polynomials} \dedicatory{Dedicated to Mark Krein on the occasion
of his 100th anniversary}
\author{D.~Alpay}
\address{Department of Mathematics \\
Ben-Gurion University of the Negev \\
P.O. Box 653 \\
84105 Beer-Sheva, Israel} \email{\tt dany@math.bgu.ac.il}
\author{A.~Dijksma}
\address{Department of Mathematics \\
University of Groningen \\
P.O. Box 407 \\
9700 AK Groningen, The Netherlands} \email{\tt
a.dijksma@math.rug.nl}

\author{H.~Langer}
\address{Institute of Analysis and Scientific Computation\\
Vienna University of Technology\\
 Wiedner Hauptstrasse 8--10 \\
A-1040 Vienna, Austria} \email{\tt
hlanger@mail.zserv.tuwien.ac.at}

\subjclass{Primary 47A57, 47B32, Secondary 42C05, 30E05}%[2000]
\keywords{Schur transformation, Nevanlinna function, realization,
symmetric operator, self-adjoint operator, moment problem,
reproducing kernel Hilbert space, orthogonal polynomials}
\thanks{D. Alpay acknowledges with thanks the Earl Katz
family for endowing the chair which supported  this research. The research of A. Dijksma and H. Langer was supported in part by the Center for Advanced Studies in Mathematics (CASM) of the Department of Mathematics,
Ben-Gurion University.} 
\begin{abstract}
A Nevanlinna function is a function which is analytic in the open upper
half plane and has a non-negative imaginary part there. In this
paper we study a fractional linear transformation for a Nevanlinna
function $n$ with a suitable asymptotic expansion at $\infty$,
that is an analogue of the Schur transformation for contractive
analytic functions in the unit disc. Applying the transformation
$p$ times we find a Nevanlinna function $n_p$ which is a
fractional linear transformation of the given function $n$. The
main results concern the effect of this transformation to the
realizations of $n$ and $n_p$, by which we mean their
representations through resolvents of self-adjoint operators in
Hilbert space. Our tools are block  operator matrix
representations, $u$--resolvent matrices, and reproducing kernel
Hilbert spaces.
\end{abstract}

\maketitle \setcounter{equation}{0}

%%%%%%%%%%%%%%%%%%%%%%%%%%%%%%%%%%%%%%%%%%%
\section{Introduction}
%%%%%%%%%%%%%%%%%%%%%%%%%%%%%%%%%%%%%%%%%%%

In the papers \cite{adls} and \cite{survey} the Schur
transformation for generalized Nevanlinna functions with a
reference point $z_1$ in the open upper half plane was considered.
An analogous transformation for Nevanlinna functions (for the
definition of a Nevanlinna function see Section \ref{Schurtrsf})
and for the reference point $\infty$ is defined in \cite[Lemma
3.3.6]{akh}, see \cite{adlline}. This transformation or a simple
modification of it we call here the {\it Schur transformation for
Nevanlinna functions}, and it is the starting point for the
present paper. To give more details, we consider a Nevanlinna
function $n$ which has for some integer $p\geq 1$ an
asymptotic expansion of order $2p+1$ at $\infty$, for example
\begin{equation}\label{i1}
n(z)=-\dfrac{s_0}{z}-\dfrac{s_1}{z^2}-\cdots-
\dfrac{s_{2p}}{z^{2p+1}}+{\rm o}\left(\dfrac
1{z^{2p+1}}\right),\quad z={\rm i}y,\ y\to\pm\infty.
\end{equation}
The Schur transform $\widehat n$ of $n$ is the function
\begin{equation}\label{i2}
\widehat n(z):=-\frac{s_0}{n(z)}-z+\frac{s_1}{s_0};
\end{equation}
the relation between $n$ and $\widehat n$ can also be written as
$$
n(z)=-\dfrac{s_0}{z-\dfrac{s_1}{s_0}+{\widehat n}(z)}.
$$
The transformed function $\widehat n=:n_1$ is again a Nevanlinna
function, but in general with an asymptotic expansion of the form
\eqref{i1} of lower order $2p-1$, and if $p>1 $ the Schur
transformation can be again applied to $ n_1$ etc. As a result we
obtain a finite sequence of Nevanlinna functions
$n_1=\widehat n,n_2=\widehat n_1,\dots,n_p=\widehat n_{p-1}$; this is the sequence of functions that
appears in the asymptotic expansion of $n$ by continued fractions,
see \cite[Section 3.3.6]{akh}.

The transformation \eqref{i2} is closely related to the finite
Hamburger moment problem. We recall that the Nevanlinna function
$n$ with an asymptotic expansion \eqref{i1} admits an integral
representation
%\begin{equation}\label{i5}
$$n(z)=\int_{-\infty}^\infty\,\dfrac{d\sigma(t)}{t-z},\quad z\ne
z^*, $$
%\end{equation}
where $\sigma$ is a bounded non-decreasing
function on $\mathbb R$. The coefficients $s_j$ in \eqref{i1} are
the {\it moments} of the function $\sigma$:
\begin{equation}\label{Achmom}
s_j=\int_{-\infty}^\infty\,t^j\,d\sigma(t),\quad j=0, 1,\dots,2p.
\end{equation}
The moment problem we have in mind is the problem to determine all
Nevanlinna functions $n$ with an expansion \eqref{i1} and given
coefficients $s_j,\,j=0,1,\dots,p$, see \cite[Theorem 3.2.1]{akh}.

An essential feature in our studies are operator representations or so-called realizations
of Nevanlinna functions, see \cite{kl1}, \cite{dls93},
 and \cite{dlls}. In fact, if the
Nevanlinna function $n$ admits an asymptotic expansion \eqref{i1}
its operator representation  takes the simple form
\begin{equation}\label{i3}
n(z)=\big((A-z)^{-1}u,u\big),\quad z\ne z^*,
\end{equation}
with some Hilbert space $\mathcal H$ with inner product
$(\cdot,\cdot)$, $u\in \mathcal H$, and a self-adjoint operator
$A$ in $\mathcal H$. We study the corresponding operator
representation of the Schur transform $\widehat n$, and also of
the functions $n_2, \dots,n_p$. For example, the function
$\widehat n$ admits an operator representation of the form
\eqref{i3} with a Hilbert space $\widehat{\mathcal H}$, an
operator $\widehat A$, and an element $\widehat u$ which are the
orthogonal complement of the element $u$ in $\mathcal H$, the
compression of $A$ to $\widehat{\mathcal H}$, and a multiple of
the projection of $Au$ onto $\widehat{\mathcal H}$, respectively.
After applying the Schur transformation $p$ times, the resulting
function $n_p$ admits an operator representation of the form
\eqref{i3} with the space
$$
\mathcal H_p'=\mathcal H \ominus\mathcal H_p,\quad \mathcal H_p:={\rm span}\,\{u,Au,\dots,A^{p-1}u\},
$$
 the operator that is
 the compression of $A$ to this space, and an element which is a
multiple of the projection of $A^pu$ onto $\mathcal H_p'$.

Since $n_p$ is obtained by   subsequent application of  fractional
linear transformations of the form \eqref{i2}, there is a
fractional linear relation between the function $n$ and the
transformed function $n_p$. We derive an explicit form for the
defining $2\times 2$ matrix function $V$ of this relation in three
ways: By calculating the resolvent of the operator $A$ in its
$2\times 2$ block matrix operator form corresponding to the
decomposition $\mathcal H=\mathcal H_p\oplus\mathcal H_p'$, by
means of the description of all generalized resolvents of a
certain symmetric restriction of $A$ with defect one in the space
$\mathcal H_{p+1}$, and via reproducing kernel methods using the
non-negative Nevanlinna kernel
$$
L_n(z,w)=\dfrac{n(z)-n(w)^*}{z-w^*},\quad z,w\in \mathbb
C\setminus\mathbb R,\ z\ne w^*.
$$

For the Nevanlinna function $n$ with an asymptotic expansion
\eqref{i1}, polynomials $e_j$ and $d_j$, $j=1,2,\dots,p$, of first
and second kind can be defined by the well-known formulas, see \cite[Chapter I]{akh}. Recall
that $e_j$ is a polynomial of degree $j$, and that $d_j$ is a
polynomial of degree $j-1$. We show  that the polynomials
$\widehat e_j$ of first kind of the transformed function $\widehat
n$ coincide, up to constant factor, with the polynomials
$d_{j+1}(z)$ of second kind for the given function $n$, whereas
the polynomials $\widehat d_j$ of second kind for $\widehat n$ are linear
combinations of $e_{j+1}$ and $d_{j+1}$. As a consequence, the
polynomials of second kind for $n$ are orthogonal with respect to
the measure generated by the non-decreasing function
$\widehat{\sigma}$ in the representation of the form \eqref{i3} of
the Nevanlinna function $\widehat n$; in this statement $\widehat
n$ can be replaced by the function $-1/n$. As in the classical
moment problem, the $2\times 2$ matrix function $V$, which
determines the fractional linear relation between $n$ and $n_p$,
can be represented by the polynomials of first and second kind.

A  short synopsis is as follows. The Schur transformation is
defined in the next section. We start with weaker forms of
the asymptotic expansion \eqref{i1}, for example
$$
n(z)=-\dfrac{s_0}{z}-\dfrac{s_1}{z^2}+{\rm o}\left(\dfrac
1{z^2}\right),\quad z={\rm i}y,\ y\to\pm\infty,
$$
and consider also a weaker form of the Schur transformation. In
Section \ref{saopreal} we mention three concrete forms of the
operator representation of $n$. The basic result of this section
is Theorem \ref{reduct} which describes the operator model for the
transformed function. Higher order approximations and the
corresponding polynomials of first and second kind are introduced
in Section \ref{orthpol}. In the operator model an asymptotic
expansion \eqref{i1} can be characterized by the fact that
$u\in{\rm dom}A^p$. The main result of this section is the
relation between the polynomials of first and second kind of $n$
and $\widehat n$ which was mentioned above. The reduction via a
$p$--dimensional subspace, that corresponds to $p$ subsequent
applications of the Schur transformation, is given in Section
\ref{reduction} by means of a block operator matrix representation
of $A$. In Section \ref{orthogpol}  the corresponding
transformation matrix $V$ is expressed in terms of the polynomials
of first and second kind. Although the final formulas are well
known (see for example \cite{akh}) this approach seems to be new.

In Section \ref{uresolvent}, applying the theory of $u$--resolvent
matrices, we derive a representation of a transformation matrix
 in an explicit form by means of the given moments; it corresponds
to Potapov's formula for the solution matrix of the Nevanlinna -
Pick problem, compare also \cite{abdl}. Finally, in Section \ref{DQO} we explain the connection between
$n$ and $n_p$ through some basic results from the theory of
resolvent invariant reproducing kernel spaces, and give another
proof for the representation of the transformation matrix
by orthogonal polynomials.

%%%%%%%%%%%%%%%%%%%%%%%%%%%%%%%%%%%%%%%%%%
\section{The Schur transformation} \label{Schurtrsf}
%%%%%%%%%%%%%%%%%%%%%%%%%%%%%%%%%%%%%%%%%

1. A Nevanlinna function is a complex  function $n$ which is
defined and analytic in the upper half plane $\mathbb C^+$ and has
the property
$$
z\in\mathbb C^+\ \  \Longrightarrow\ \  {\rm Im}\,\ n(z)\ge 0.
$$
We always suppose that $n$ is extended to the lower half plane
$\mathbb C^-$ by the relation
\begin{equation}\label{symm}
n(z)=n(z^*)^*, \quad z\in\mathbb C^-,
\end{equation}
 and to those points of the real axis into which it can be continued analytically. The set of all Nevanlinna functions is denoted by $\mathbf N_0$. Recall that $n\in \mathbf N_0$ if and only if $n$ is analytic in $\mathbb C\setminus \mathbb R$ and the kernel
$$
L_n(z,w)=\dfrac{n(z)-n(w)^*}{z-w^*},\quad z,w\in \mathbb C\setminus \mathbb R,\ z\ne w^*,
$$
is positive definite.

Let $n \in \mathbf N_0$ and  consider the following
properties of $n$:
\begin{enumerate}
\item[$(1_0)$]$ n(z)=-\dfrac{s_0}{z}+{\rm o}\left(\dfrac
{1}{z}\right)$,\vspace*{2mm} \item[$(2_0)$]$
 n(z)=-\dfrac{s_0}z+{\rm O}\left(\dfrac 1{z^2}\right),
$\vspace*{2mm}
\item[$(3_0)$]$
 n(z)=-\dfrac{s_0}z-\dfrac{s_1}{z^2}+{\rm o}\left(\dfrac
 1{z^2}\right),
$
\end{enumerate}
where  here and in the
following, the limit relations are understood to hold for
$z\to\pm\infty$ along the imaginary axis. The assumption \eqref{symm} implies that  $s_0$ and $s_1$ are real numbers. Evidently,
$(3_0)\Longrightarrow (2_0)\Longrightarrow (1_0)$. The function
$n$ satisfies the assumption $(1_0)$ if and only if it belongs to
the class $(R_0)$ of \cite{KK}, which means that it admits an
integral representation
\begin{equation}\label{int}
n(z)=\int_{-\infty}^{+\infty}\frac 1 {t-z}\, d\sigma(t),\quad
z\in\mathbb C\setminus\mathbb R,
\end{equation}
where $\sigma$ is a bounded non-decreasing function on $\mathbb
R$. Then
$$
\int_{-\infty}^{+\infty}d\sigma(t)=s_0,
$$
hence $s_0\geq 0$, and if $s_0=0$ then $n(z)\equiv 0$. With the
representation \eqref{int} of $n$ the assumption
\begin{equation}\label{card}
\int_{-\infty}^{\infty}|t|\,d\sigma(t)<\infty
\end{equation}
implies that $(3_0)$ is satisfied. Indeed, \eqref{card} implies
that
%\begin{equation}\label{c11}
$$
s_1=\int_{-\infty}^\infty t\, d\sigma(t)
$$
%\end{equation}
exists and with $z=iy$
$$
z^2\left(n(z)+\dfrac{s_0}{z}+\dfrac{s_1}{z^2}\right)
=\int_{-\infty}^\infty \dfrac{t^2}{t-z}\,d\sigma(t)=
\int_{-\infty}^\infty \dfrac{t^2+iyt}{t^2+y^2}\,t\,d\sigma(t) ={\rm
o}(1).
$$

The assumptions $(1_0), (2_0)$, and $(3_0)$ are all different. To
see that $(1_0)\not \Longrightarrow (2_0)$ we show that if $n \in\mathbf N_0$ has the representation \eqref{int} with
$ {\rm supp} \,\sigma\subset [0,\infty)$ and
$$
 \displaystyle\int_0^\infty
d\sigma(t)<\infty,\quad\displaystyle\int_0^\infty\, t\,
d\sigma(t)=\infty,
$$
 then $(2_0)$ does not hold: Let $c>0$ be given
arbitrarily (large) and choose $K>0$ such that $\displaystyle\int_0^K
t\,d\sigma(t)\ge c$. If $y$ is chosen large enough then for $0\le
t\le K$ we have
$$
\dfrac{y^2}{t^2+y^2}\ge \frac12,
$$
and hence
$$
\displaystyle\int_{-\infty}^\infty \frac {
y^2t}{t^2+y^2}\,d\sigma(t)\ge \frac c2,
$$
 and therefore, with $z=iy$,
$$
z^2\left(n(z)+\dfrac{s_0}{z}\right)=
z\int_{-\infty}^\infty \frac t {t-z}\,d\sigma(t)
=\int_{-\infty}^\infty \frac
{-y^2t+iyt^2}{t^2+y^2}\,d\sigma(t)\neq {\rm O}(1),
$$
which implies
that $(2_0)$ does not hold. Thus, for example, the function
$$
n(z)=\dfrac{-1}{z-\sqrt{-z}}=\int_0^\infty
\dfrac{1}{t-z}\,\dfrac{dt}{\pi (t+1)\sqrt{t}}
$$
satisfies $(1_0)$ but not $(2_0)$.

Let $n$ be the  Nevanlinna function, defined in the upper half plane by
$$
n(z)=\dfrac{-s_0}{z+\gamma+f(z)}, \ z \in \mathbb C^+,
$$
where $s_0$ is a positive real number, $\gamma$ is a complex
number with ${\rm Im}\,\gamma>0$, and $f$ is a Nevanlinna function
such that $f(z)={\rm o}(1)$. It has the properties
$$
\lim_{z=iy, \, y \rightarrow \infty}\,
z^2\left(n(z)+\dfrac{s_0}{z}\right)=\gamma s_0, \qquad
\lim_{z=iy,\, y \rightarrow -\infty}\,
z^2\left(n(z)+\dfrac{s_0}{z}\right)=\gamma^* s_0,
$$
and hence $n$ satisfies $(2_0)$ but, since the two limits are
different (and non-real), it  does not satisfy $(3_0)$.

Instead of
the assumption $(3_0)$ also the assumption
\begin{equation}\label{ca}
n(z)=-\dfrac{s_0}z-\dfrac{s_1}{z^2}+{\rm O}\left(\frac
1{z^3}\right)
\end{equation}
seems reasonable. However, according to \cite[Bemerkung
1.11]{kl1}, \eqref{ca} implies the existence of a real number
$s_2$ such that
\begin{equation}\label{bemerkung}
n(z)=-\frac{s_0}z-\frac{s_1}{z^2}-\frac{s_2}{z^3}+{\rm o}
\left(\frac 1{z^3}\right);
\end{equation}
this relation will be considered in Section \ref{orthpol} as assumption $(1_1)$. The implication \eqref{ca} $\Longrightarrow$
\eqref{bemerkung} can also be seen from the integral
representation \eqref{int} of $n$: \eqref{ca} implies
$$
z^3\left(\int_{-\infty}^\infty
\dfrac{1}{t-z}+\dfrac{1}{z}\,d\sigma(t)+ \dfrac{s_1}{z^2}\right)
={\rm O}(1),$$ and hence with $z=iy$
$$-\int_{-\infty}^\infty \dfrac{y^2t(t+iy)}{t^2+y^2}
\,d\sigma(t)+iys_1={\rm O}(1), \quad y \rightarrow \infty.$$
Taking the imaginary part we see that
$$
s_1={\rm lim}_{y \rightarrow \infty}\,\int_{-\infty}^\infty
\dfrac{y^2t}{t^2+y^2} \,d\sigma(t)
$$
and taking the real part we see that there exist real numbers $C$
and $y_0$ such that
$$\int_{-\infty}^\infty \dfrac{y^2t^2}{t^2+y^2}
\,d\sigma(t)\leq C, \quad y\geq y_0.
$$
This implies that
$$s_2:=\int_{-\infty}^\infty t^2\,d\sigma(t) < \infty,$$
hence $$\int_{-\infty}^\infty |t|\,d\sigma(t) < \infty$$ and
$$s_1=\int_{-\infty}^\infty t\,d\sigma(t).$$
Now \eqref{bemerkung} easily follows from the integral
representations of $n$ and the expressions for the real numbers $s_0$, $s_1$, and
$s_2$: With $z=iy$ we have
$$
z^3\left(n(z)+\dfrac{s_0}{z}+\dfrac{s_1}{z^2}+\dfrac{s_2}{z^3}\right)=
\int_{-\infty}^\infty \dfrac{t^3}{t-z}\, d\sigma(t)=
\int_{-\infty}^\infty \dfrac{t^4+it^3y}{t^2+y^2}\, d\sigma(t)={\rm
o}(1).
$$

2. Now we define the basic transformations considered this paper.
 \begin{definition}\label{strans}
If $n\in\mathbb N_0$  satisfies the assumption $(1_0)$ or
$(2_0)$, the {\it Schur type transform $\widetilde n$ of $n$} is the function
\begin{equation}\label{ee1}
\widetilde n(z)=\frac{-s_0}{n(z)}-z,
\end{equation}
if $n\in\mathbb N_0$  satisfies the assumption $(3_0)$  the {\it
Schur transform $\widehat n$ of $n$} is the function
\begin{equation}\label{ee2}
\widehat n(z)=\frac{-s_0}{n(z)}-z+\frac{s_1}{s_0}.
\end{equation}
\end{definition}
The difference between the formulas \eqref{ee1} and \eqref{ee2} is
just in the additive real constant $s_1/s_0$: under the stronger
assumption $(3_0)$ this constant assures that the transform tends
to zero if $z$ tends to $\pm \infty$ along the imaginary axis, see
\eqref{eel3} below.

 The relations \eqref{ee1} and \eqref{ee2} can also be written as a
 first step of a continued fraction expansion
 $$
 n(z)=-\frac{s_0}{z+{\widetilde n}(z)}
,\qquad {\rm or}\ \ n(z)=-\dfrac{s_0}{z-\frac{s_1}{s_0}+{\widehat
n}(z)}.
$$
\begin{theorem}\label{el}
The following equivalences hold:
\begin{eqnarray}
n\in \mathbf N_0\ {\rm and\  satisfies}\
(1_0)&\Longleftrightarrow&
\widetilde n \in \mathbf N_0,\,\widetilde n(z)={\rm o}(z),\label{eel1}\\
\label{eel2}n\in \mathbf N_0\ {\rm and\  satisfies}\
(2_0)&\Longleftrightarrow& \widetilde n \in \mathbf
N_0,\,\widetilde n(z)={\rm O}(1),\\\label{eel3} n\in \mathbf N_0\
{\rm and\ satisfies}\ (3_0)&\Longleftrightarrow& \widehat n \in
\mathbf N_0,\,\widehat n(z)={\rm o}(1).
\end{eqnarray}
\end{theorem}

\begin{proof}
We have
$$
\widetilde n(z)=\widehat
n(z)-\dfrac{s_1}{s_0}=-\dfrac{s_0+zn(z)}{n(z)}.
$$
A straightforward calculation yields
$$
{\rm Im}\, \widetilde n(z)={\rm Im}\, \widehat n(z)=\dfrac{{\rm
Im}\, z}{|n(z)|^2}\left(s_0\dfrac{{\rm Im}\, n(z)}{{\rm Im}\,
z}-|n(z)|^2\right),
$$
and the estimate
$$
|n(z)|^2=\left|\int_{-\infty}^{+\infty} \dfrac{d\sigma(t)}{t-z}\right|^2\le
\int_{-\infty}^{+\infty}\dfrac{d\sigma(t)}{|t-z|^2}\int_{-\infty}^{+\infty}d\sigma(t)=
\dfrac{{\rm Im}\, n(z)}{{\rm Im}\, z}s_0
$$
implies $\widetilde n,\,\widehat n\in \mathbf N_0$. The asymptotic
properties of $\widetilde n$ follow in case \eqref{eel1} from the
relation
$$
\dfrac{\widetilde n(z)}z=-\dfrac{s_0}{zn(z)}-1={\rm o}(1),
$$
in case \eqref{eel2} from the relation
$$
\widetilde n(z)=-\dfrac{s_0+zn(z)}{n(z)}=\dfrac{z\,{\rm
O}\left(\dfrac 1z\right)}{zn(z)},
$$
and for $\widehat n(z)$ in case \eqref{eel3} in a similar way or
from \cite[Lemma 3.3.6]{akh}.

Conversely, starting from $\widetilde n(z)$ as in \eqref{eel1},
the relation
$$
z\left(n(z)+\dfrac {s_0}z\right)=s_0\dfrac{\dfrac{\widetilde
n(z)}{z}}{1+\dfrac{\widetilde n(z)}{z}}
$$
 implies that from $\widetilde n(z)={\rm o}(z)$ it follows that
$n$ satisfies $(1_0)$. The corresponding proofs for \eqref{eel2}
and \eqref{eel3} are similar.
\end{proof}

\section{Self-adjoint operator representations} \label{saopreal}
A function $n\in\mathbf N_0$ admits a {\it self-adjoint operator
representation} or {\it realization}
with a self-adjoint relation
$A$ in some Hilbert space $\calH$ of the form
\begin{equation}\label{itn}
n(z)=n(z_0)^* + (z-z_0^*)\left(\left(I+(z-z_0)(A-z)^{-1}\right)
v,v\right)
\end{equation}
with $z_0$ an arbitrary non-real  number $z_0$ and an element $v\in
\calH$, see \cite{kl1}, \cite{dls93},
 and \cite{dlls}. If $v$ is chosen to be a
 {\it generating element for} $A$, which means that
 $$
 \calH={\overline{\rm span}}\,\left\{v, (A-z)^{-1}v\,\big|\,z\in\mathbb
C\setminus\mathbb R\right\}
 $$
 and which is always possible, then the operator representation
 \eqref{res} is called {\it minimal} and then it is unique up to
 unitary equivalence. We have the following equivalences, see \cite{lt}:
\begin{eqnarray*}
 A \text{ is an operator} &\Longleftrightarrow& n(z)={\rm o}(z),\\
v \in {\rm dom} A & \Longleftrightarrow& {\rm lim}_{y \rightarrow \infty}y\,{\rm Im}\,n(iy)<\infty;
\end{eqnarray*}
 for $n\in\mathbf N_0$ the latter limit always exists: it is either a 
non-negative number or
 $\infty$.

If the Nevanlinna function $n$ satisfies the assumption $(1_0)$
(or any of the assumptions $(2_0),\,(3_0)$) the representation
\eqref{itn} can be simplified to
\begin{equation}\label{res}
n(z)=\big((A-z)^{-1}u,u\big),\qquad z\in\mathbb C\setminus\mathbb R,
\end{equation}
where $A$ is a self-adjoint operator in some Hilbert space
 $\calH,\,u\in \calH,\,(u,u)=s_0$. If $u$ is chosen to be a
 generating element for $A$, or equivalently,
 $$
 \calH={\overline{\rm span}}\,\left\{(A-z)^{-1}u\,\big|\,z\in\mathbb
C\setminus\mathbb R\right\}
 $$
 which is always possible, then the operator representation
 \eqref{res} is also called {\it minimal} and then it is unique up to
 unitary equivalence. The representation \eqref{res} follows from
 \eqref{itn} and the above mentioned equivalences by taking
 $u=c(A-z_0)v$ with some unimodular complex number $c$.

Here are three examples for a more concrete choice of the triplet
$\mathcal H,\,A,\,u$ in  \eqref{res} for the given function
$n\in\mathbf N_0$ with integral representation \eqref{int}.

\begin{itemize}

\item[(1)] $\mathcal H=L^2(\sigma)$, $A$ is the operator of
multiplication with the independent variable, and $u(t)\equiv
1,\,t\in\mathbb R$.

\item[(2)] $\calH$ is the completion  of the  linear span of the functions ${\bf r}_z,\,z\in
\mathbb C \setminus \mathbb R$:
$$
{  r}_z(t):=\frac{1}{t-z},\quad t\in\mathbb R,
$$
 with inner product defined by
$$
({  r}_z,{  r}_\zeta)=\dfrac{n(z)-n(\zeta)^*}{z-\zeta^*},\quad
z,\zeta\in \mathbb C \setminus \mathbb R,\ z\ne\zeta^*,
$$
 $A$ is the operator of multiplication by $t$, and
$u(t)\equiv 1,\,t\in\mathbb R$.

\item[(3)]  $\mathcal H$ is the reproducing kernel Hilbert space
$\mathcal L(n)$ with reproducing kernel
$$
L_n(z,w)=\dfrac{n(z)-n(w)^*}{z-w^*},\quad z,w\in \mathbb C \setminus \mathbb R,\ z\ne w^*,
$$
$A$ is the self-adjoint operator whose resolvent $(A-z)^{-1}$ is
the difference-quotient operator $R_z$:
$$
(R_zf)(\zeta)=\dfrac{f(\zeta)-f(z)}{\zeta-z}, \quad f \in \mathcal
L (n),
$$
and take $u=n$; this function belongs to the space $\mathcal
L(n)$, since $n$ satisfies the condition $(1_0)$. Recall that the reproducing
property of the kernel $L_n$ is reflected in the inner product of the
space ${\mathcal L}(n)$:
\[
\langle f,L_n(\,\cdot\,, z\rangle_{{\mathcal L}(n)}=f(z),\quad f\in{\mathcal L}
(n),\ z\in{\mathbb C}\setminus{\mathbb R}.
\]
That \eqref{res}
holds follows from
$$(R_zn)(\zeta)=L_n(\zeta, z^*)$$
and the reproducing property of the kernel $L_n$:
\begin{eqnarray*}
\lefteqn{\langle(A-z)^{-1}u,u\rangle_{\mathcal L(n)}=\langle
R_zn,n\rangle_{\mathcal L(n)}}\\&&\hspace{2cm}  = \langle
n,L_n(\,\cdot\,,z^*)\rangle_{\mathcal L(n)}^*=n(z^*)^*=n(z).
\end{eqnarray*}
\end{itemize}

The unitary equivalence  of  the representations in (1) and (2)
follows easily from the relation
\[
({  r}_z,{
r}_\zeta)=\int_{-\infty}^{\infty}\,\dfrac{d\sigma(t)}{(t-z)(t-\zeta^*)},\quad
z,\zeta\in\mathbb C\setminus\mathbb R,
\]
and the fact that the functions ${  r}_z,\,z\ne z^*$, form a total
set in $L^2(\sigma)$. The unitary equivalence between the two representations of $n$ in
(2) and (3) is given by the mapping $U$:
$$
 U( r_z)= L_n(\,\cdot, z^*);
 $$
in particular, we have $Uu=n$ where $u$ is the function
$u(t)\equiv 1,\,t\in\mathbb R$. The space $L^2(\sigma)$
(or the equivalent space in (2)) we denote also by $\calH(n)$. We mention that the definition of the spaces in (2) and (3) can also be used for
generalized Nevanlinna functions, whereas in this case the space $L^2(\sigma)$
need not be defined. In Sections \ref{orthpol} -  \ref{uresolvent} we will
prove theorems using the representation of $n$ in (1), in Section \ref{DQO} we reprove some of these
results using the representation in the reproducing kernel Hilbert
space $\mathcal L(n)$.

Since, according to Theorem \ref{el}, the functions $\widetilde n$ and $\widehat n$ in Definition
\ref{strans} belong to the class
$\mathbf N_0$ and are ${\rm o}(z)$ for $z=iy,\,y\to\infty$, they admit again an operator
representation of the form \eqref{itn}, for example,
\begin{equation}\label{itntilde}
\widehat n(z)=\widehat n(z_0)^* + (z-z_0^*)\left((\widehat
A-z_0)(\widehat A-z)^{-1}\widehat v,\widehat v\right)
\end{equation}
with a self-adjoint operator $\widehat A$ in some Hilbert space
$\widehat{\mathcal H}$, $z_0$ an arbitrary non-real  number, and
an element $\widehat v \in \widehat{\mathcal H}$. Clearly, as the
difference between the functions $\widetilde n$ and $\widehat n$
is just an additive real constant, the operator representation for
$\widetilde n$ can be chosen the same, that is, in
\eqref{itntilde} $\widehat n$ can be replaced by $\widetilde n$.

\begin{theorem} \label{reduct}
Let $n\in \mathbf N_0$ satisfying the condition $(1_0)$ and with operator representation \eqref{res} be given, and let
$$
%\begin{equation}\label{ee2}
\widehat n(z)=\frac{-s_0}{n(z)}-z+\frac{s_1}{s_0}
%\end{equation}
$$
 be the Schur transform of $n$ from \eqref{ee2}. Then in the operator
representation \eqref{itntilde} of $\widehat n$ we can choose
$\widehat\calH=\{u\}^\perp,\,\widehat A$ in $\widehat \calH$ as
the compression of $A$ to $\widehat{\mathcal H}: \widehat
A=\widehat PA|_{\widehat\calH}$, where $\widehat P$ is the
orthogonal projection in $\mathcal H$ onto $\widehat{\mathcal H}$,
and the element $\widehat v$ as
$$
\widehat v=\dfrac{\|u\|}{((A-z_0)^{-1}u,u)}\widehat
P(A-z_0)^{-1}u.
$$
%If $n$ satisfies even the condition $(1_1)$ then
If $\widehat n$ also satisfies the condition $(1_0)$
\footnote{This is the case when $n$ satisfies condition $(1_1)$
defined in Section \ref{orthpol}, see Lemma \ref{asymplemma}.},
then
$$ \widehat n(z)=\left((\widehat A - z)^{-1}\widehat u,\widehat u)
 \right), \quad
 \widehat u:=\dfrac{\widehat PAu}{\|u\|}.
 $$
\end{theorem}

\begin{remark}
The resolvent of $\widehat A$ is given by
$$
(\widehat A-z)^{-1}=(A-z)^{-1}-
\dfrac{((A-z)^{-1}\,\cdot\,,u)}{((A-z)^{-1}u,u)}(A-z)^{-1}u,
$$
and
$$
\widehat v=\dfrac{\|u\|(A-z_0)^{-1}u-((A-z_0)^{-1}u,u)\dfrac
u{\|u\|}}{|((A-z_0)^{-1}u,u)|}.
$$
Note that $\left((A-z_0)^{-1}u,u\right)=n(z_0)\neq 0$, otherwise
$n(z)\equiv0$.
\end{remark}

\begin{proof}[Proof of Theorem {\rm \ref{reduct}}]
(1) Suppose  that $n$ satisfies $(1_0)$. Then we have
$\|u\|=\sqrt{s_0}$ and
\begin{equation}\label{uuu}
\widehat n(z)-\widehat
n(z_0)^*=-\dfrac{\|u\|^2}{r(z)}-z+\frac{\|u\|^2}{r(z_0)^*}+z_0^*,
\end{equation}
where we have put $r(z):=((A-z)^{-1}u,u)$. It remains to show that the  expression on the right hand side of \eqref{uuu} equals
$$
(z-z_0^*)\left((\widehat A-z_0)(\widehat A-z)^{-1})\widehat
v,\widehat v\right)=(z-z_0^*)\|\widehat
v\|^2+(z-z_0^*)(z-z_0)\left((\widehat A-z)^{-1}\widehat v,\widehat
v\right).
$$
This is a straightforward calculation, we only indicate some formulas:
\begin{eqnarray*}
(z-z_0)(\widehat A-z)^{-1}\widehat v&=&\frac{\|u\|}{|r(z_0)|}\left(\frac{r(z_0)}{r(z)}(A-z)^{-1}u-(A-z_0)^{-1}u\right),\\
\|\widehat
v\|^2&=&\|u\|^2\frac{\|(A-z_0)^{-1}u\|^2}{|r(z_0)|^2}-1,
\end{eqnarray*}
and  \begin{eqnarray*} \lefteqn{(z-z_0^*)(z-z_0)\left((\widehat
A-z)^{-1}\widehat v,\widehat v\right)}\\&& =\|u\|^2\left(\frac 1
{r(z_0)^*}-\frac1{r(z)}-\frac{(z-z_0^*)\|(A-z_0)^{-1}u\|^2}{|r(z_0)|^2}\right).
\end{eqnarray*}
(2) Now assume that $\widehat n$  satisfies $(1_0)$. Then $\widehat v
\in {\rm dom}\,\widehat A \subset {\rm dom}\,A$ and the equality
$$\widehat v=\dfrac{\|u\|}{r(z_0)}(A-z_0)^{-1}u-\dfrac{u}{\|u\|}$$
shows that also $u \in {\rm dom}\,A$. If we take $\widehat
u=-(\widehat A-z_0)\widehat v$, then (see after \eqref{res})
$\widehat n$ has the asserted representation. It remains to show
that $\widehat u=\widehat PAu/\|u\|$. We have
$$
\widehat u=-(A-z_0)\widehat v=\dfrac{1}{\|u\|}(A-z_0)u-
\dfrac{\|u\|}{r(z_0)}u.
$$
Taking the inner product of both sides with $u$ and using
$(\widehat u,u)=0$, we see that
$$\dfrac{\|u\|}{r(z_0)}=\dfrac{\left((A-z_0)u,u\right)}{\|u\|^3}$$
and hence
$$\widehat
u=\dfrac{1}{\|u\|}(A-z_0)u-\dfrac{\left((A-z_0)u,u\right)}{\|u\|^3}u
=\dfrac{1}{\|u\|}\left(Au-\dfrac{(Au,u)}{\|u\|^2}u\right)
=\dfrac{1}{\|u\|}\widehat P Au.$$

\end{proof}

\section{Higher order asymptotics. Orthogonal polynomials}\label{orthpol}

1. For $n\in
\mathbf N_0$ and some integer $p \geq 1$ we introduce the assumptions
\begin{enumerate}
\item[$(1_p)$]$
n(z)=-\dfrac{s_0}{z}-\dfrac{s_1}{z^2}-\cdots-\dfrac{s_{2p}}{z^{2p+1}}+{\rm
o}\left(\dfrac 1{z^{2p+1}}\right)$,\vspace*{2mm}
\item[$(2_p)$]$n(z)=-\dfrac{s_0}{z}-\dfrac{s_1}{z^2}-\cdots-\dfrac{s_{2p}}{z^{2p+1}}
+{\rm O}\left(\dfrac 1{z^{2p+2}}\right)$,\vspace*{2mm}
\item[$(3_p)$]$n(z)=-\dfrac{s_0}{z}-\dfrac{s_1}{z^2}-\cdots-
\dfrac{s_{2p+1}}{z^{2p+2}}+{\rm o}\left(\dfrac
1{z^{2p+2}}\right)$.
\end{enumerate}
Again, $(3_p)\,\Longrightarrow\, (2_p)\, \Longrightarrow\, (1_p)$,
and  by \cite[Satz 1.10]{kl1} for the operator representation
the assumption $(1_p)$ is equivalent to $u\in{\rm dom} \,A^p$. That
is, for the above representation with the space $\calH(n)$  the functions
$$
{\bf t}_k(t):=t^k,\quad k=0,1,\dots,p,
$$
belong to $\calH(n)$ and the first $p$ of these elements, ${\bf
t}_0,{\bf t}_1,\dots,{\bf t}_{p-1}$, belong to ${\rm dom}\,A$.
Moreover, the formal relation
$$
(A-z)^{-1}=-\sum_{j=0}^{2p}\frac{A^j}{z^{j+1}}+\frac{A^{2p+1}}{z^{2p+1}}(A-z)^{-1}
$$
implies easily
\begin{eqnarray*}
n(z)&=&\left((A-z)^{-1}u,u\right)\\&=&-
\sum_{j=0}^{p}\frac{(A^ju,u)}{z^{j+1}}-\sum_{j=p+1}^{2p}
\frac{(A^nu,A^{j-p}u)}{z^{j+1}}+\frac1{z^{2p+1}}\left(A(A-z)^{-1}A^pu,A^pu\right).
\end{eqnarray*}
It follows that
$$
s_j=\left\{\begin{array}{lll}
(A^ju,u)&\ {\rm if}&\ j=0,1,\dots,p,\\[2mm]
(A^nu,A^{j-p}u)&\ {\rm if}&\ j=p+1,p+2,\dots,2p.
\end{array}\right.
$$
Therefore the above assumptions are equivalent to the following relations for the operator $A$ and the generating element $u,\,u\in{\rm dom}\,A^p$:
\begin{equation}\label{ess}
\begin{array}{rcl}
(1_p)\!\!\!&\Longleftrightarrow&\!\!\! u\in{\rm dom}\,A^p,
\\[2mm]
 (2_p)\!\!\!& \Longleftrightarrow& \!\!\!u\in{\rm dom}\,A^p,\
 z(A(A-z)^{-1}A^pu,A^pu)={\rm O}(1),\\[2mm]
 (3_p)\!\!\!&\Longleftrightarrow&\!\!\!u\in{\rm dom}\,A^p,\
 z(A(A-z)^{-1}A^pu,A^pu)+\alpha={\rm o}(1)\ {\rm with\ }\alpha\in\mathbb R;
\end{array}
\end{equation}
in fact, in the last equivalence we have $\alpha=s_{2p+1}$.

Now we consider a function $n\in \mathbf N_0$ with the property $(1_p)$ for some $p>1$.
 For $0\le k\le p$, by $S_k$ we denote the
 $(k+1)\times(k+1)$ Hankel matrix
\begin{equation}\label{hankel}
S_k:=
\begin{pmatrix} s_0&s_1&\cdots&s_k\\
s_1&s_2&\cdots&s_{k+1}\\
\vdots&\vdots& &\vdots\\
s_{k}&s_{k+1}&\cdots&s_{2k}
\end{pmatrix};
\end{equation}
it is the Gram matrix associated with the $k+1$ functions $\mathbf
t_0, \mathbf t_1, \ldots, \mathbf t_{k}$, and we introduce the Gram
determinants
\begin{equation}\label{dk}
D_{k}:={\rm det}\, S_k= \left| \begin{array}{cccc} s_0&s_1&\cdots&s_{k}\\
s_1&s_2&\cdots&s_{k+1}\\
\vdots&\vdots& &\vdots\\
s_{k}&s_{k+1}&\cdots&s_{2k}
\end{array}\right|,\quad k=0,1,\dots p.
\end{equation}
Further, for $k=1,\dots,p$, $\mathcal H_k$ denotes the
$k$--dimensional subspace
$$
\mathcal H_k:={\rm span}\, \{{\bf t}_0,{\bf t}_1,\dots,{\bf
t}_{k-1}\}
$$
of $\calH (n)$. Evidently, the subspace $\mathcal H_k$ is
non-degenerated if and only if $D_{k-1}\ne 0$.

In the rest of this section we suppose
that $D_{p-1}\ne 0$, that is, the subspace $\mathcal H_p$ is
non-degenerated. If $D_p=0$, then the function $n$ with the given
asymptotics $(j_p)$ is uniquely determined and rational of Mac
Millan degree $p$, in fact, see \cite[pp. 22,23]{akh}
$$
n(z)=- \dfrac{d_p(z)}{e_p(z)},
$$
where the polynomials $e_p$ of degree $p$ and $d_p$ of degree
$p-1$ are defined below. To exclude this (simple) case we often suppose that even $D_p\ne
0$; clearly, this implies $D_{p-1}\ne 0$.

 As a basis in $\mathcal H_p$ we choose a
system of elements $e_k\in\calH(n)=L^2(\sigma),\,k=0,1,\dots,p-1$, which is
obtained from the system ${\bf t}_0,{\bf t}_1,\dots,{\bf t}_p$ by
the Gram--Schmidt orthonormalization procedure. This so-called system of
{\it orthogonal polynomials of the
first kind, associated with the function $n$} is defined by the following properties, $j,k=0,1,\dots,p-1$:
\begin{enumerate}
\item $e_0(z)\equiv \displaystyle 1/{\sqrt{s_0}}$,\vspace*{1mm}
\item $e_k(z)$ is a real polynomial of degree $k$ with positive
leading coefficient,\vspace*{1mm} \item $ (e_j,e_k)= \delta_{jk}.
$
\end{enumerate}
Then,  see \cite[(1.4)]{akh},
\begin{equation}\label{innere}
e_k(z)= \frac 1{\sqrt{D_{k-1}D_k}}\left| \begin{array}{cccc}
s_0&s_1&\cdots&s_k\\
s_1&s_2&\cdots&s_{k+1}\\
\vdots&\vdots&\cdots&\vdots\\
s_{k-1}&s_k&\cdots&s_{2k-1}\\
1&z&\cdots&z^k \end{array} \right|,\quad k=1,2,\dots,p-1,
\end{equation}
and by this formula with $ k=p$ also a polynomial $e_p$ can be defined.
Evidently, ${e}_p\in\mathcal H_p^{[\perp]}$, and
$$
{\rm span}\,\{{\bf t}_0,{\bf t}_1,\dots,{\bf t}_k\}={\rm span}\,
\{e_0,\,e_1,\,\dots,\,e_k\},\quad k=0,1,\dots,p.
$$

The orthogonal polynomials $e_j,\,j=0,1,\dots,p$, satisfy the difference equations
\begin{equation}\label{pol1}
b_{k-1}e_{k-1}(z)+a_ke_k(z)+b_ke_{k+1}(z)=ze_k(z),\quad k=0,1,\dots,p-1,
\end{equation}
with real numbers $a_k,\,k=0,1,\dots,p-1,\,b_{-1}=0$, and positive
numbers $b_k,\,k=1,\dots,p-1$, and the `initial condition' $e_0(z)=
\displaystyle 1/{\sqrt{s_0}}$. Explicit formulas for $a_k,\,b_k$
can be given, see \cite{akh}; we note that
\begin{equation} \label{a0b0}
 a_0=\dfrac{s_1}{s_0}, \quad b_0=\dfrac{\sqrt{s_2s_0-s_1^2}}{s_0}.
 \end{equation}

The relation \eqref{pol1}
implies that with respect to the basis $e_0,e_1,\dots,e_{p-1}$ of
the space $\calH_p$ the compression $A_p$ of the operator $A$ to
$\calH_p$ is given by the Jacobi matrix
\begin{equation}\label{ja}
\calA_p:=\left(
\begin{array}{cccccc}
a_0&b_0&0&\cdots&0&0\\
b_0&a_1&b_1&\cdots&0&0\\
\vdots&\vdots&\vdots&\vdots&\vdots&\vdots\\
0&0&0&\cdots&a_{p-2}&b_{p-2}\\0&0&0&\cdots&b_{p-2}&a_{p-1}
\end{array}\right),
\end{equation}
and that
$$
Ae_{p-1}=b_{p-2}e_{p-2}+a_{p-1}e_{p-1}+b_{p-1}e_p.
$$
The latter relation means for the orthogonal polynomials
$$
b_{p-2}e_{p-2}(z)+a_{p-1}e_{p-1}(z)+b_{p-1}e_p(z)=ze_{p-1}(z),
$$
therefore, the eigenvalues of $\calA_p$ are the zeros of the
polynomial $e_p$. For later use we write the last $p-1$ difference
equations \eqref{pol1} explicitly in the form
\begin{equation}\label{op1}
\left\{\begin{array}{lcl}
%a_0e_0+b_0e_1&=&ze_0\\
b_0e_0+a_1e_1+b_1e_2&=&ze_1\\
b_1e_1+a_2e_2+b_2e_3&=&ze_2\\
b_2e_2+a_3e_3+b_3e_4&=&ze_3\\
&\vdots&\\
b_{p-2}e_{p-2}+a_{p-1}e_{p-1}+b_{p-1}e_p&=&ze_{p-1};
\end{array}\right.
\end{equation}
this system of homogeneous equations for $e_0,e_1,\dots,e_p$
determines the orthogonal polynomials uniquely if we add the
initial conditions
\begin{equation}\label{ic1}
e_0(z)=\dfrac1 {\sqrt{s_0}},\quad
e_1(z)=\dfrac{z-a_0}{b_0\sqrt{s_0}};
\end{equation}
the second condition is just the first equation in \eqref{pol1}.

The {\it polynomials of the second kind, associated with the
function $n\in\mathbf N_0$}, are the functions
$d_k,\,k=0,1,\dots,p$, defined as follows:
\begin{equation}\label{innerd}
d_k(z)=\sqrt{s_0}\left(\displaystyle\frac{e_k(z)-e_k(\cdot)}{z-\cdot},e_0\right)=
\left(\displaystyle\frac{e_k(z)-e_k(\cdot)}{z-\cdot},u\right),\quad
k=0,1,\dots,p.
\end{equation}
Hence $d_0(z)=0$ and $d_k$ is a polynomial of degree $k-1$, $k\geq
1$. The definition of $d_k$ and the relation \eqref{pol1} imply
that
\begin{equation}\label{pol2}
b_{k-1}d_{k-1}(z)+a_kd_k(z)+b_kd_{k+1}(z)=zd_k(z),\quad
k=1,\dots,p-1.
\end{equation}
Therefore the polynomials $e_k$
and $d_k$
 satisfy for $k=1,2,\dots,p-1$ the same difference equations  but with
 different  initial conditions:
 \begin{equation}\label{ic}
  \quad d_0(z)=0,\quad  d_1(z)= \dfrac{\sqrt{s_0}}{b_0}.
 \end{equation}
 For later use we write the difference equations \eqref{pol2} in the form
 \begin{equation}\label{op2}
\left\{\begin{array}{lcl}
a_1d_1+b_1d_2&=&zd_1\\
b_1d_1+a_2d_2+b_2d_3&=&zd_2\\
b_2d_2+a_3d_3+b_3d_4&=&zd_3\\
&\vdots&\\
b_{p-2}d_{p-2}+a_{p-1}d_{p-1}+b_{p-1}d_p&=&zd_{p-1}.\\
\end{array} \right.
\end{equation}

For any two solutions $u_0,\dots,u_p$ and $v_0,\dots,v_p$ of the
difference equations \eqref{pol1} with $b_{-1}=0$:
$$
\begin{array}{rcl}
zu_k(z)&=&b_{k-1}u_{k-1}(z)+a_ku_k(z)+b_ku_{k+1}(z),\\[2mm]
\zeta
v_k(\zeta)&=&b_{k-1}v_{k-1}(\zeta)+a_kv_k(\zeta)+b_kv_{k+1}(\zeta),
\end{array}
\quad k=0,1,\dots,p-1,
$$
the Christoffel--Darboux formulas hold:
\begin{equation}\label{sc}
\begin{array}{rcl}
\displaystyle\sum_{k=m}^{p-1}(z-\zeta)u_k(z)v_k(\zeta)&=&b_{p-1}\big(u_p(z)v_{p-1}
(\zeta)-u_{p-1}(z)v_p(\zeta)\big)\\
&&-b_{m-1}\big(u_m(z)v_{m-1} (\zeta)-u_{m-1}(z)v_m(\zeta)\big);
\end{array}
\end{equation}
in particular,
\begin{equation}\label{lundy}
d_p(z)e_{p-1}(z)-e_p(z)d_{p-1}(z)=\frac 1{b_{p-1}}.
\end{equation}

2. In this subsection we assume that $n\in\mathbf N_0$ satisfies the
assumption $(1_p)$ for some $p\ge 1$, and we consider its Schur transform $\widehat
n$ from \eqref{ee2}. For the following lemma see \cite[Lemma
2.1]{der2}, we sketch the proof.

\begin{lemma} \label{asymplemma}
Suppose that $n$ satisfies $(1_p)$ for some $p\ge 1:$
$$
n(z)=-\dfrac{s_0}{z}-\dfrac{s_1}{z^2}-\cdots-
\dfrac{s_{2p}}{z^{2p+1}}+{\rm o}\left(\dfrac{1}{z^{2p+1}}\right),
$$
then its Schur transform $\widehat n$ satisfies $(1_{p-1}):$
$$
%\cn(z):=\
%-\dfrac 1{n(z)}-\dfrac{z}{s_0}+\dfrac{s_1}{s_0^2}\
\widehat n(z)=\ -\dfrac{\widehat s_0}{z}-\dfrac{\widehat
s_1}{z^2}-\cdots-\dfrac{\widehat s_{2p-2}}{z^{2p-1}}+{\rm
o}\left(\dfrac1{z^{2p-1}}\right),
$$
with for $j=0,1,\ldots, 2p-2$
\begin{equation}\label{coeffnhat}
\widehat s_j=\dfrac{(-1)^{j+1}}{s_0^{j+2}}\left|
\begin{array}{cccccc}
s_1&s_0&0&\cdots&0&0\\
s_2&s_1&s_0&\cdots&0&0\\
\vdots&\vdots&\ddots&\ddots&\vdots&\vdots\\
s_{j+1}&s_{j}&s_{j-1}&\cdots&s_1&s_0\\
s_{j+2}&s_{j+1}&s_j&\cdots&s_2&s_1
\end{array}\right|,
\end{equation}

\end{lemma}

\begin{proof}
Write
$$
\widehat n(z)=z\left(1+\dfrac{s_1}{s_0z}+\dfrac{
s_2}{s_0z^2}+\cdots+\dfrac{s_{2p}}{s_0z^{2p}}+{\rm
o}\left(\dfrac1{z^{2p}}\right)\right)^{-1}-z+\dfrac{s_1}{s_0}.
$$
If we set
$$
q(z)=\dfrac{s_1}{s_0z}+\dfrac{
s_2}{s_0z^2}+\cdots+\dfrac{s_{2p}}{s_0z^{2p}}+{\rm
o}\left(\dfrac1{z^{2p}}\right),
$$
then
\begin{eqnarray*}\lefteqn{\dfrac{1}{1+q(z)}=1-q(z)+ \cdots
-q(z)^{2p-1}+\dfrac{q(z)^{2p}}{1+q(z)}}\\&& \hspace{1cm} =1-q(z)+
\cdots -q(z)^{2p-1} +{\rm o}\left(\dfrac{1}{z^{2p-1}}\right)
\end{eqnarray*}
and
\begin{eqnarray*}
\widehat n(z)&=& \dfrac z {1+q(z)}-\left(z-
\dfrac{s_1}{s_0}\right)\dfrac{1+q(z)}{1+q(z)}\\
&=&\left(-zq(z)+\dfrac{s_1}{s_0}(1+q(z))\right)
\dfrac 1{1+q(z)}\\
&=&\left(-z\left(\dfrac{ s_2}{s_0z^2}+\cdots+ \dfrac{
s_{2p}}{s_0z^{2p}}+{\rm o}\left(\dfrac 1{z^{2p}}\right)\right)
 + \dfrac{ s_1}{s_0}q(z)\right)\dfrac 1{1+q(z)}\\
&=&\left(-\dfrac{ s_2}{s_0z}-\cdots- \dfrac{
s_{2p}}{s_0z^{2p-1}}+{\rm o} \left(\dfrac
1{z^{2p-1}}\right)+\dfrac{s_1}{s_0}q(z)\right) \dfrac 1{1+q(z)},
\end{eqnarray*}
which is of the needed form. Formula \eqref{coeffnhat} for the
coefficients $\widehat s_j$ can be obtained by equating powers of
$z$ from both sides of the equality
$$n(z)\left(z-\dfrac{s_1}{s_0}+\widehat n(z)\right)=-s_0.$$
\end{proof}

 Now we can formulate the main result of this
subsection.

\begin{theorem}\label{orthogpolnhat}
Let $n\in\mathbf N_0$ satisfy condition $(1_p)$ for some $p\ge 2$.
%, and consider the function
%$$
%\cn(z)=-\dfrac 1{n(z)}-\dfrac z{s_0} +\dfrac{s_1}{s_0^2}.
%$$
If $e_k$ and $d_k,\,k=0,1,\dots,p$, denote the polynomials of
first and second kind associated with the function $n$, and
$\widehat e_k$ and $\widehat d_k,\,k=0,1,\dots,p-1$, denote the
polynomials of first and second kind associated with the Schur
transform $\widehat n$ of $n$, then for $k=0,1, \dots,p-1$ the
following relations hold:
%\begin{eqnarray*}
\begin{equation}\label{rel1}
\widehat e_k(z)=\dfrac{1}{\sqrt{s_0}}d_{k+1}(z),
\end{equation}
\begin{equation}\label{rel2}
\begin{array}{rcl}
\widehat
d_k(z)\!\!\!&=&\!\!\!\!-\sqrt{s_0}e_{k+1}(z)+\dfrac{1}{\sqrt{s_0}}
\left(z-\dfrac{s_1}{s_0}\right)d_{k+1}(z)\\[4mm]
\!\!\!&=&\!\!\!\!b_0\left(\dfrac{e_{k+1}(z)-e_{k+1} (\cdot)}{z-\cdot}, e_1
\right)\!=\!\dfrac{1}{\sqrt{s_0}}\left(\dfrac{e_{k+1}(z)-e_{k+1}(\cdot)}{z-\cdot},
\cdot -a_0\right).
\end{array}
\end{equation}
\end{theorem}
\begin{remark} (i) Here we write $\widehat e$ and $\widehat d$ for the
polynomials of first and second kind associated with the Schur
transform $\widehat n$ of $n$, but the reader is reminded that
these functions are {\it not} the Schur transforms of the
polynomials $e$ and $d$. \\[1mm]
(ii) If $\check e_k$ and $\check d_k$ stand for the polynomials of
first and second kind associated with the Nevanlinna function
$-1/n$, then for $k=0,1, \dots,p-1$
$$\check e_k(z)=d_{k+1}(z), \quad \check
d_k(z)=-e_{k+1}(z)+\dfrac{1}{s_0}
\left(z-\dfrac{s_1}{s_0}\right)d_{k+1}(z).
$$
The first equality readily follows from the fact that the spectral
functions of $\widehat n$ and $-1/n$ only differ by a factor
$s_0$. The second equality can be obtained by tracing the proof
below; the only difference lies in \eqref{sat}: with evident
notation, it should be replaced by
$$
\check s_0=\dfrac{s_0s_2-s_1^2}{s_0^3}.
$$
\end{remark}

\begin{proof}[Proof of Theorem {\rm \ref{orthogpolnhat}}]
For the function $\widehat n$, again with evident notation,  we
have
\begin{equation}\label{sat}
\widehat s_0=\dfrac{s_0s_2-s_1^2}{s_0^2}
\end{equation}
 and, as a consequence of Theorem \ref{reduct},
\begin{equation}\label{AA}
\widehat{\mathcal A}_{p-1}\!=\!
\begin{pmatrix}
\widehat a_0\!&\!\widehat b_0\!&\!0\!&\!\cdots\!\!&\!0\!\!&\!0\\
\widehat b_0&\widehat a_1&\widehat b_1&\cdots&0&0\\
\vdots&\vdots&\vdots&\vdots&\vdots&\vdots\\
0&0&0&\cdots&\widehat a_{p-3}&\widehat b_{p-3}\\
0&0&0&\cdots&\widehat b_{p-3}&\widehat a_{p-2}
\end{pmatrix}
\!=\!
\begin{pmatrix}
\!a_1\!&\!b_1\!&\!0\!&\!\cdots\!&\!0\!&\!0\!\\
\!b_1\!&\!a_2\!&\!b_2\!&\!\cdots\!&\!0&\!0\!\\
\!\vdots\!&\!\vdots\!&\!\vdots\!&\!\vdots\!&\!\vdots\!&\!\vdots\!\\
\!0&\!0&\!0\!&\!\cdots\!&\!a_{p-2}\!&\!b_{p-2}\!\\\!0&\!0&\!0\!&\!\cdots\!&\!b_{p-2}\!&\!a_{p-1}\!
\end{pmatrix}.
\end{equation}
For the  $\widehat e_j,\,j=0,\dots,p-1,$ we find
$$
\widehat e_0(z)=\dfrac 1{\sqrt{\widehat
s_0}}=\dfrac{s_0}{\sqrt{s_0s_2-s_1^2}}=\dfrac{1}{b_0}=
\dfrac{1}{\sqrt{s_0}}\,d_1(z)
$$
and $\widehat e_1,\widehat e_2,\dots,\widehat e_{p-1}$ follow from
the equations (see \eqref{pol1}):
%\begin{equation}\label{hop2}
$$
\left\{\begin{array}{lcl}
a_1\widehat e_0+b_1\widehat e_1&=&z\widehat e_0\\
b_1\widehat e_0+a_2\widehat e_1+b_2\widehat e_2&=&z\widehat e_1\\
b_2\widehat e_1+a_3\widehat e_2+b_3\widehat e_3&=&z\widehat e_2\\
&\vdots&\\
b_{p-2}\widehat e_{p-3}+a_{p-1}\widehat e_{p-2}+b_{p-1}\widehat
e_{p-1} &=&z\widehat e_{p-2}.
\end{array}\right.
$$
%\end{equation}
Since these equation coincide with \eqref{op2} we obtain
$$
\widehat e_j=c\,d_{j+1},\quad j=0,1,\dots,p-1.
$$
The constant $c$ can be determined from the initial condition
$\widehat e_0=c\,d_1$, which gives $c=1/\sqrt{s_0}$.
Therefore
$$
\widehat e_j(z)=\dfrac{1}{\sqrt{s_0}}\,d_{j+1}(z),\quad
j=0,1,\dots,p-1,
$$
and \eqref{rel1} is proved.

For the polynomials of second kind
 $\widehat d_j$ we obtain in a similar way
$$
\widehat d_0(z)=0,\quad \widehat d_1(z)=\dfrac{\sqrt{\widehat
s_0}}{\widehat b_0},
$$
%\begin{equation}\label{hop2}
$$
\left\{\begin{array}{lcl} \widehat a_1\widehat d_1+\widehat
b_1\widehat d_2&=
&z\widehat d_1\\
\widehat b_1\widehat d_1+\widehat a_2\widehat d_2+
\widehat b_2\widehat d_3&=&z\widehat d_2\\
\widehat b_2\widehat d_2+\widehat a_3\widehat d_3+
\widehat b_3\widehat d_4&=&z\widehat d_3\\
&\vdots&\\
\widehat b_{p-3}\widehat d_{p-3}+\widehat a_{p-2} \widehat
d_{p-2}+\widehat b_{p-2}\widehat d_{p-1}&=
&z\widehat d_{p-2}\\
\end{array}\right.
$$
%\end{equation}
(one equation less than in \eqref{op2}). These
equations can be written as
\begin{equation}\label{hop22}
\left\{\begin{array}{lcl}
a_2\widehat d_1 + b_2\widehat d_2&=&z\widehat d_1 \\
 b_2\widehat d_1+ a_3\widehat d_2+ b_3\widehat d_3&
 =&z\widehat d_2\\
 b_3\widehat d_2+ a_4\widehat d_3+ b_4\widehat d_4&
 =&z\widehat d_3\\
&\vdots&\\
 b_{p-2}\widehat d_{p-3}+ a_{p-1}\widehat d_{p-2}+
 b_{p-1}\widehat d_{p-1}&=&z\widehat d_{p-2}.
 \end{array} \right.
\end{equation}
The last $p-3$ equations of this system coincide with the last
$p-3$  equations of \eqref{op1} and \eqref{op2}. Therefore a
solution vector $(\widehat d_j)_1^{p-1}$ of the last $p-3$
equations of \eqref{hop22} can be obtained as a linear combination
of the solution vectors $(e_j)_2^p$ and $(d_j)_2^p$ of the last
$p-3$ equations in \eqref{op1} and \eqref{op2}:
$$
\widehat d_j=\gamma e_{j+1} + \delta d_{j+1},\quad
j=1,2,\dots,p-1.
$$
Now $\gamma,\,\delta$ have to be found such that these relations
hold also for $j=0$ with $\widehat d_0(z)=0$, and for $j=1$ with
$\widehat d_1(z)=\dfrac{\sqrt{\widehat s_0}}{\widehat b_0}.$ Since $\widehat d_0(z)=0$ it follows that
$$
0=\gamma e_1(z) + \delta d_1(z)=\gamma\,\dfrac{z-a_0}{b_0} \dfrac
1{\sqrt{s_0}}+\delta\,\dfrac{\sqrt{s_0}}{b_0}=
\gamma\,\dfrac{z-\dfrac{s_1}{s_0}}{b_0}\dfrac 1{\sqrt{s_0}}
+\delta\,\dfrac{\sqrt{s_0}}{b_0},
$$
which is satisfied for
$$
\gamma = -\varepsilon,\quad\delta=\varepsilon
\left(\dfrac z{s_0}-\dfrac{s_1}{s_0^2}\right).
$$
The relation $\widehat d_1(z)= \dfrac{\sqrt{\widehat
s_0}}{\widehat b_0}$ implies
\begin{eqnarray*}
\widehat d_1(z)&=&\varepsilon\left(-e_2(z) + \left(\dfrac z{s_0} -
\dfrac{a_0}{s_0}\right)d_2(z)\right)\\
&=&\varepsilon\left(-\dfrac{(z-a_1)e_1 - b_0e_0}{b_1} +
\dfrac {z-a_0}{s_0} \dfrac{(z-a_1)d_1}{b_1}\right)\\
&=&\varepsilon\left(-\dfrac{(z-a_1)\dfrac{z-a_0}{b_0}
\dfrac 1{\sqrt{s_0}} - b_0\dfrac 1{\sqrt{s_0}}}{b_1} +
\dfrac{z-a_0}{s_0} \dfrac{(z-a_1)
\dfrac{\sqrt{s_0}}{b_0}}{b_1}\right)\\
&=&\varepsilon\,\dfrac{b_0}{b_1\sqrt{s_0}} =\dfrac{\sqrt{\widehat
s_0}}{\widehat b_0}=\dfrac{b_0}{\widehat b_0},
\end{eqnarray*}
hence
$$
\varepsilon =\dfrac{b_1}{\widehat b_0}\sqrt{s_0}.
$$
According to \eqref{AA} we find $\varepsilon =\sqrt{s_0}$. This
proves the first equality in \eqref{rel2}. The remaining
equalities follow from \eqref{innerd} and the second equality in
\eqref{ic1}.
\end{proof}

3. In this subsection we give a second proof of Theorem
\ref{orthogpolnhat} using asymptotic expansions, see \cite[(1.34b)]{akh}: Assume that $n \in
\mathbf N_0$ satisfies $(1_p)$ for some $p\geq 2$, that is,
$$
n(z)=-\dfrac{s_0}{z}-\dfrac{s_1}{z^2}-\cdots-\dfrac{s_{2p}}{z^{2p+1}}+{\rm
o}\left(\dfrac 1{z^{2p+1}}\right),
$$
then
\begin{equation}\label{134b}
-\dfrac{d_p(z)}{e_p(z)}=-\dfrac{s_0}{z}-\dfrac{s_1}{z^2}-
\cdots-\dfrac{s_{2p-1}}{z^{2p}}+{\rm O}\left(\dfrac
{1}{z^{2p+1}}\right).
\end{equation}
According to \cite[the second to last formula on p. 22]{akh} the
function $-d_p/e_p$ is a Nevanlinna function and by
\cite[Bemerkung 1.11]{kl1} there is a real number $t_{2p}$ such
that $-d_p/e_p$ has the asymptotic expansion
\begin{equation}\label{134bKL}
-\dfrac{d_p(z)}{e_p(z)}=-\dfrac{s_0}{z}-\dfrac{s_1}{z^2}-
\cdots-\dfrac{s_{2p-1}}{z^{2p}}-\dfrac{t_{2p}}{z^{2p+1}}+ {\rm
o}\left(\dfrac 1{z^{2p+1}}\right).
\end{equation}
By Lemma \ref{asymplemma},
$$
\widehat n(z)= -\dfrac{\widehat s_0}{z}-\dfrac{\widehat
s_1}{z^2}-\cdots-\dfrac{\widehat s_{2p-2}}{z^{2p-1}}+{\rm
o}\left(\dfrac1{z^{2p-1}}\right),
$$
and hence by \eqref{134b}
\begin{equation}\label{134bcheck}
-\dfrac{\widehat d_{p-1}(z)}{\widehat e_{p-1}(z)}=-\dfrac{\widehat
s_0}{z}-\dfrac{\widehat s_1}{z^2}- \cdots-\dfrac{\widehat
s_{2p-3}}{z^{2p-2}}+{\rm O}\left(\dfrac {1}{z^{2p-1}}\right).
\end{equation}
By Lemma \ref{asymplemma}, the Schur transform of the function
$-d_p/e_p$ in \eqref{134bKL} has the asymptotic expansion
\begin{eqnarray}\nonumber\lefteqn{
\widehat{\left(-\dfrac{d_p}{e_p}\right)}(z)=
\dfrac{s_0e_p(z)}{d_p(z)} -z
+\dfrac{s_1}{s_0}=:\dfrac{r(z)}{d_p(z)}}\\&& \label{last=}
=-\dfrac{\widehat s_0}{z}-\dfrac{\widehat
s_1}{z^2}-\cdots-\dfrac{\widehat
s_{2p-3}}{z^{2p-2}}-\dfrac{\widehat t_{2p-2}}{z^{2p-1}} +{\rm
o}\left(\dfrac1{z^{2p-1}}\right),\\&& = -\dfrac{\widehat
s_0}{z}-\dfrac{\widehat s_1}{z^2}- \cdots-\dfrac{\widehat
s_{2p-3}}{z^{2p-2}}+{\rm O}\left(\dfrac {1}{z^{2p-1}}\right),
\nonumber
\end{eqnarray}
where only the number $\widehat t_{2p-2}$ depends on $t_{2p}$
according to formula \eqref{coeffnhat}. Here the polynomial $r$,
defined via the second equality sign, is given by
$$r(z)=s_0e_p(z)-\left(z-\dfrac{s_1}{s_0}\right)d_p(z)$$
and its degree is $\leq p$. Comparing \eqref{134bcheck} with
\eqref{last=}, we find that
$$
\dfrac{r(z)}{d_p(z)}-\dfrac{\widehat d_{p-1}(z)}{\widehat
e_{p-1}(z)}= {\rm O}\left(\dfrac {1}{z^{2p-1}}\right).
$$
The degree of the product $d_p\,\widehat e_{p-1}$ equals $2p-2$
and hence
$$
\dfrac{r(z)}{d_p(z)}=\dfrac{\widehat d_{p-1}(z)}{\widehat
e_{p-1}(z)},
$$
which readily implies that for some number $k \neq 0$
\begin{equation}\label{check}
\widehat e_{p-1}(z)=kd_p(z), \quad \widehat
d_{p-1}(z)=k\left(s_0e_p(z)-\left(z-\dfrac{s_1}{s_0}\right)d_p(z)\right).
\end{equation}
We claim $k=1/\sqrt{s_0}$. With the proof of the claim the proof
of the theorem is complete.

To prove the claim we note that the
leading coefficient of the polynomial $e_k$ is equal to
$\sqrt{D_{k-1}/D_k}$ and that, by \eqref{pol1},
$$\sqrt{\dfrac{D_k}{D_{k+1}}}=\dfrac{1}{b_k}
\sqrt{\dfrac{D_{k-1}}{D_k}}.$$ Hence
$$\sqrt{\dfrac{D_{p-1}}{D_p}}=\dfrac{1}{b_{p-1}}\cdots
\dfrac{1}{b_1}\sqrt{\dfrac{D_0}{D_1}}=\dfrac{1}{b_{p-1}}\cdots
\dfrac{1}{b_1}\dfrac{1}{b_0}\dfrac{1}{\sqrt{s_0}},
$$
and, similarly, because of \eqref{AA} and \eqref{sat},
$$
\sqrt{\dfrac{\widehat D_{p-2}}{\widehat D_{p-1}}}
=\dfrac{1}{\widehat b_{p-2}}\cdots \dfrac{1}{\widehat
b_0}\dfrac{1}{\sqrt{\widehat s_0}}=\dfrac{1}{b_{p-1}}\cdots
\dfrac{1}{b_1}\dfrac{1}{b_0}.
$$ From \eqref{check} we obtain
$$
\dfrac{1}{k}\,\sqrt{\dfrac{\widehat D_{p-2}}{\widehat D_{p-1}}}= \lim_{z \rightarrow \infty}\dfrac{d_p(z)}{z^{p-1}}= s_0\lim_{z \rightarrow \infty}\dfrac{e_p(z)}{z^p}= s_0\,\sqrt{\dfrac{D_{p-1}}{D_p}},
$$
that is,
$$
\sqrt{s_0}\,\dfrac{1}{b_{p-1}}\cdots \dfrac{1}{b_1}\dfrac{1}{b_0}=
\dfrac{1}{k}\dfrac{1}{b_{p-1}}\cdots
\dfrac{1}{b_1}\dfrac{1}{b_0}.
$$
Therefore, $k=1/\sqrt{s_0}$ and the claim holds.

 \section{Reduction via a $p$-dimensional subspace}
\label{reduction}
 Let again $n\in\mathbf N_0$ with the property
$(1_p)$ be given. We decompose the space $\mathcal H(n)$ with $\mathcal H_p={\rm span}\, \{{\bf t}_0,{\bf t}_1,\dots,{\bf
t}_{p-1}\}$ as
follows:
\begin{equation}\label{uu}
 \mathcal H(n)=\mathcal H_p \oplus\mathcal H'_p.
\end{equation}
Then, evidently, $e_p\in\mathcal H'_p$. The corresponding matrix
representation of the operator $A$ is
\begin{equation}\label{pp}
A=\left(\begin{array}{cc} A_0&\widetilde B\\B&D
\end{array}\right)
\end{equation}
with $A_0$  given by the Jacobi matrix $\calA_0$ from \eqref{ja},
$$
B=b_{p-1}(\,\cdot\,,e_{p-1}){e}_{p},\quad\widetilde B=
b_{p-1}(\,\cdot\,,{e}_p)e_{p-1}.
$$
 The operator $A_0$ is bounded and
self-adjoint in $\calH_p$, and $D$ is self-adjoint (possibly
unbounded) in $\calH'_p$. In the following theorem we express the
function $n$ by means of the entries of the matrix in \eqref{pp}.
We set
$$
a_{00}(z):=\big((A_0-z)^{-1}u,u\big),\quad
a_{11}(z):=((A_0-z)^{-1}e_{p-1},e_{p-1}),
$$
\begin{equation}\label{ioi}
a(z):=((A_0-z)^{-1}u,e_{p-1})=((A_0-z)^{-1}e_{p-1},u),
\end{equation}
and
$$
\check r(z):=\begin{vmatrix}\, a_{00}(z) & a(z) \\
\,a(z) & a_{11}(z) \end{vmatrix}.
$$
The last equality in \eqref{ioi} follows from $a(z^*)^*=a(z)$, in fact, by
Cramer's rule, $a(z)=\sqrt{s_0}/(a_0-z)$ if $p=1$ and
$a(z)=(-1)^{p-1}\sqrt{s_0}\,b_0\ldots b_{p-2}/{\rm det}\,(\mathcal
A_0-z)$ if $p \geq 2$.

\begin{theorem}\label{jk+p}
Suppose that the function $n\in \mathbf N_0$ satisfies for some
$p\ge 1$ one of the assumptions $(j_p), j=1,2,3$, and that $D_p\ne
0$, see \eqref{dk}.
 Then
\begin{equation}\label{le}
n(z)=\frac{\check r(z) n_p(z)-a_{00}(z)} {a_{11}(z)n_p(z)-1},
\end{equation}
where
$n_p(z):=((D-z)^{-1} u_p, u_p)$ with $u_p:= b_{p-1}e_p$. The function $n_p$ belongs to
$\mathbf N_0$ and satisfies the assumption $(j_0)$.
Moreover, for $k\ge 1$
we have
\begin{equation}\label{rad}
u\in {\rm dom}\,A^{p+k} \Longleftrightarrow u_p\in {\rm
dom}\,D^{k},
\end{equation}
and $n$ satisfies the assumption $(j_{p+k})$ if and only if $n_p$
satisfies the assumption $(j_k)$.
\end{theorem}

\begin{remark} If the operator representation \eqref{res} of $n$ is $n(z)=\big((A-z)^{-1}u,u\big)$ with the space $\mathcal H=\mathcal H(n)$,  then, according to Theorem \ref{reduct},  the operator representation of $n_1=\widehat n$ is given by
$$
\mathcal H_1'=\mathcal H \ominus\{u\},\quad A_1=P_{\mathcal
H_1'}A|_{\mathcal H_1'},\quad u_1=\dfrac{P_{\mathcal
H_1'}Au}{\|u\|},
$$
 the operator representation of $n_2=\widehat n_1$ by
$$
\mathcal H_2'=\mathcal H \ominus\{u, Au\},\quad A_2=P_{\mathcal
H_2'}A|_{\mathcal H_2'},\quad u_2=\dfrac{P_{\mathcal
H_2'}A^2u}{\|u\|\|u_1\|},
$$
 and, via induction, the operator representation of $n_p=\widehat n_{p-1}$ by
 $$
 \mathcal H_p'=\mathcal H \ominus\{u,Au, \cdots,A^{p-1}u\},\quad
 A_p=P_{\mathcal H_p'}A|_{\mathcal H_p'},\quad u_p=\dfrac{P_{\mathcal H_p'}A^pu}{\|u\|\|u_1\|\ldots \|u_{p-1}\|}.
 $$
 Note that $\|u\|=\sqrt{s_0}$ and, by Theorem \ref{jk+p},
 $\|u_j\|=b_{j-1}$, $j=1, \ldots,p$.
\end{remark}

\begin{proof}[Proof of Theorem {\rm \ref{jk+p}}]
With the matrix \eqref{pp}, the equation $(A-z)x=u$ becomes
\begin{eqnarray}\label{cc}
(A_0-z)x_1\,+\,b_{p-1}(x_2, e_p)e_{p-1}&=&u,\\
b_{p-1}(x_1,e_{p-1}) e_p\,+\ \  (D-z)x_2\ &=&0,\nonumber
\end{eqnarray}
where $x$ is written as $x=\begin{pmatrix} x_1 & x_2 \end{pmatrix}^\top$ according to the decomposition \eqref{uu}. The second
equation implies
$$
x_2=-b_{p-1}(x_1,e_{p-1})(D-z)^{-1} e_p.
$$
We insert this into  \eqref{cc} and apply $(A_0-z)^{-1}$ to get
\begin{equation}\label{hh}
x_1=(A_0-z)^{-1}u+b_{p-1}^2(x_1,e_{p-1})\left((D-z)^{-1} e_p,
e_p\right)(A_0-z)^{-1}e_{p-1}.
\end{equation}
Now take the inner product with $e_{p-1}$ and solve the obtained
equation for $(x_1,e_{p-1})$:
$$
(x_1,e_{p-1})=\frac{\left((A_0-z)^{-1}u,e_{p-1}\right)}{1-b_{p-1}^2
\left((D-z)^{-1}e_p, e_p\right)
\left((A_0-z)^{-1}e_{p-1},e_{p-1}\right)}.
$$
Observing that  $ n(z)=\left((A-z)^{-1}u,u\right)=(x_1,u)$,  the
relation \eqref{hh} yields \eqref{le}. %In the matrix \eqref{pp}
%the entries with possible exception of $D$ are bounded operators.
To prove \eqref{rad}, denote by $P'$ the orthogonal projection onto $\calH'_p$ in
\eqref{uu}. If $k=1$, then $u\in{\rm  dom}\,A^{p+1}$ is equivalent to
$v:=A^{p}u\in{\rm dom}\,A$. Since $D$ is the only entry in the matrix of \eqref{pp} which is possibly unbounded ($B$ and $\widetilde B$ are even one--dimensional), $v\in A$ is equivalent to the fact that the non-zero
component  $P'v$, which is a multiple of $e_p$, belongs to ${\rm
dom}\, D$. If $k=2$ we observe that
\begin{equation}\label{qu}
A^2=\begin{pmatrix}
A_0^2+\widetilde BB&A_0\widetilde B+\widetilde BD\\
BA_0+DB&B\widetilde B+D^2
\end{pmatrix}.
\end{equation}
 Since $e_p\in{\rm dom}\,D$ the operators $DB$ and $\widetilde BD$ and
 hence all the entries in the matrix representation of $A^2$
 except possibly $D^2$ are bounded.
 Now $u\in{\rm  dom}\,A^{p+2}$ is equivalent to
$v=A^{p}u\in{\rm dom}\,A^2$, and hence, by \eqref{qu}, $P'v\in{\rm
dom}\,D^2$. The claim for arbitrary $k$ follows by induction.

The last claim of the theorem for $j=1$ follows immediately from
\eqref{rad} and the first equivalence in \eqref{ess}. For $j=2,3$ we also use the
equivalences in \eqref{ess}.  A simple calculation yields
\begin{equation}\label{respp}
(A-z)^{-1}=\begin{pmatrix}
R_{11}(z)&R_{12}(z)\\R_{12}(z^*)^*&R_{22}(z)
\end{pmatrix}
\end{equation}
with
$$
R_{11}(z)=S_1(z)^{-1},\quad
R_{12}(z)=-b_{p-1}\left((D-z)^{-1}\,\cdot\,,e_p\right)S_1(z)^{-1}e_{p-1},
$$
$$
R_{22}(z)=(D-z)^{-1}+ b_{p-1}^2
\left(S_1(z)^{-1}e_{p-1},e_{p-1}\right)\left((D-z)^{-1}\,\cdot\,,e_p\right)(D-z)^{-1}e_p,
$$
where
$S_1(z):=A_0-z-b_{p-1}^2\left((D-z)^{-1}e_p,e_p\right)(\,\cdot\,,e_{p-1})e_{p-1}$,
the first Schur complement. It is easy to see that for $f,g\in\calH_p$ we have
$$
\lim_{y\to\infty}\,{\rm i} y\big(S_1({\rm i} y)^{-1}f,g\big)=-(f,g).
$$
Now we observe the relation
$$
A(A-z)^{-1}=
\begin{pmatrix}
A_0R_{11}(z)+\widetilde B R_{12}(z^*)^*&A_0R_{12}(z)+\widetilde BR_{22}(z)\\[2mm]
BR_{11}(z)+DR_{12}(z^*)^*&BR_{12}(z)+DR_{22}(z)
\end{pmatrix}
$$
and the fact that for $z={\rm i}y,\,y\to\infty$, for example
$zA_0R_{11}(z)=zA_0S_1(z)$ has a limit,
\begin{eqnarray*}
zA_0R_{12}(z)&\!=\!&-zb_{p-1}\big((D-z)^{-1}\cdot,e_p\big)A_0S_1(z)^{-1}e_{p-1}={\rm o}(1),\\
zDR_{22}(z)&\!=\! &zD(D\!-\!z)^{-1}+zb_{p-1}^2\big(S_1(z)^{-1}e_p,e_p\big)\big((D\!-\!z)^{-1}\cdot,e_p\big)D(D\!-\!z)^{-1}e_p\\
&\!=\!&zD(D-z)^{-1}+{\rm o}(1),
\end{eqnarray*}
etc. These relations imply for example with $v=A^{p+k}u$
$$
z\left(A(A-z)^{-1}v,v\right)=z\left(D(D-z)^{-1}P'v,P'v\right)+{\rm
O}(1).
$$
Since $P'v\in{\rm span}\{e_p,e_{p+1},\dots,e_{p+k}\},\,e_p,e_{p+1},\dots,e_{p+k-1}\in{\rm dom}\,D$ and hence
$$
z\left(D(D-z)^{-1}x',x'\right)={\rm O}(1) \text{ for } x'\in {\rm span}\{e_p,e_{p+1},\dots,e_{p+k}\},
$$
and since $P'v$  has a non-zero component in the direction of
$e_{p+k}$ the claim follows from \eqref{ess}.
\end{proof}

\section{Representation of the
transformation matrix by orthogonal polynomials}\label{orthogpol}

The $2 \times 2$ matrix function which generates the fractional
linear transformation \eqref{le} we denote in the following by
$V$:
\begin{equation} \label{v}
V(z):=\dfrac{1}{a(z)} \begin{pmatrix}
\check r(z)&-a_{00}(z)\\[2mm]
a_{11}(z)&-1
\end{pmatrix}
\end{equation}

In this section we express $V$ by the polynomials of first and
second kind. To this end, the elements of $\mathcal H_p$ are
considered as column vectors with respect to the basis
$e_0,e_1,\dots,e_{p-1}$.

First we solve the equation $(A_0-z)x=e_{p-1}$ in $\mathcal H_p$.
With  the Jacobi matrix $\mathcal A_0$ from \eqref{ja} this
equation becomes
$$
 \calA_0
x-zx=e_{p-1},
$$
or
%\begin{equation}\label{aa}
$$  (\mathcal A_0-z)\left(
  \begin{array}{c}\xi_1\\ \vdots\\ \xi_{p-1}\\ \xi_{p}
  \end{array} \right)=\left(
  \begin{array}{c}0\\  \vdots \\0\\1
  \end{array} \right).
$$
%\end{equation}
According to the definition of the orthogonal polynomials of first
kind, the solution of the system with the $1$ in the
 last component of the vector on the right hand side replaced by
 $-b_{p-1} e_p(z)$
 is the vector with components
 $e_0(z),\,e_1(z),\dots,e_{p-1}(z)$. It follows that
 %\begin{equation}\label{ab}
$$
 \xi_j=-\displaystyle \frac{e_{j-1}(z)}{b_{p-1}
  e_p(z)},\quad j=1,2,\dots,p,
$$
 %\end{equation}
  and hence
$$
 \begin{array}{ccccc}
   \left((A_0-z)^{-1}e_{p-1},e_{p-1}\right)\!\!\!\!&\!\!=\!\!&\!\!\!\!\left(
   x,e_{p-1}\right)\!\!\!\!
   &\!\!\!\!=\!\!\!\!&\!\!
   -\dfrac{e_{p-1}(z)}{b_{p-1}
    e_p(z)},\\[4mm]
      \left((A_0-z)^{-1}e_{p-1},e_0\right)&=&(x,e_0)&=&
      -\displaystyle\frac{e_0(z)}{b_{p-1}
 e_p(z)},
\end{array}
$$
that is,
\begin{equation}\label{hf}
a_{11}(z)=-\dfrac{e_{p-1}(z)}{b_{p-1}
    e_p(z)}, \quad
    a(z)=-\displaystyle\frac{1}{b_{p-1}e_p(z)}.
\end{equation} Next we solve the equation
$(A_0-z)x=u=\sqrt{s_0}\,e_0$. As above, in matrix form it becomes

%\begin{equation}\label{01}
$$
  (\mathcal A_0-z)\left(
\begin{array}{c}\xi_1\\ \xi_2\\ \vdots\\ \xi_{p}
\ \end{array} \right)=
\left(\begin{array}{c}\sqrt{s_0}\\0\\\vdots\\0\end{array}\right).
$$
%\end{equation}
 According to the definition of the
polynomials of the second kind and because of
$b_0\,d_1(z)=\sqrt{s_0}$ we have
$$
(\mathcal A_0-z)\left( \begin{array}{c}
0\\d_1(z)\\\vdots\\d_{p-2}(z)\\d_{p-1}(z)
\end{array}
\right)= \left( \begin{array}{c} \sqrt{s_0}\\0\\\vdots\\0\\-
b_{p-1}d_p(z)
\end{array}
\right)=\left( \begin{array}{c} \sqrt{s_0}\\0\\\vdots\\0\\
0
\end{array}
\right)-b_{p-1}d_p(z)\left( \begin{array}{c} 0\\0\\\vdots\\0\\1

\end{array}
\right).
$$
It follows that
%$$
%(\mathcal A_0-z)\left\{\left( \begin{array}{c}
%0\\d_{1}(z)\\\vdots\\d_{n-1}(z)
%\end{array}
%\right)-\frac{ d_n(z)}{ e_n(z)}\left(
%\begin{array}{c}
%e_0(z)\\e_1(z)\\\vdots \\e_{n-1}(z)
%\end{array}
%\right)\right\}=\left( \begin{array}{c} b_0d_1(z)\\0\\
%\vdots\\0
%\end{array}
%\right),
%$$
%and the solution of (\ref{01}) becomes
$$
\left(
\begin{array}{c}\xi_1\\ \xi_2\\ \vdots\\ \xi_{p}
\ \end{array} \right)=(\mathcal A_0-z)^{-1}\left( \begin{array}{c}
\sqrt{s_0}\\0\\\vdots\\0
\end{array}
\right)=\left(
\begin{array}{c} 0\\d_1(z)\\\vdots\\d_{p-1}(z)
\end{array}
\right)-\frac{ d_p(z)}{ e_p(z)}\left(
\begin{array}{c}
e_0(z)\\e_1(z)\\\vdots\\e_{p-1}(z)
\end{array}
\right)
$$
and hence
\begin{equation}\label{hf1}
a_{00}(z)=\left((A_0-z)^{-1}u,u\right)=-\frac{ d_p(z)}{ e_p(z)}.
\end{equation}
Inserting the expressions from (\ref{hf}) and (\ref{hf1}) into
(\ref{le}) and observing the relation \eqref{lundy} we find that
$V$ can be written as
%\begin{equation}\label{V}
$$
V(z)=\begin{pmatrix}
-d_{p-1}(z)&-b_{p-1}d_p(z)\\
e_{p-1}(z) & b_{p-1}e_p(z)
\end{pmatrix},
$$
%\end{equation}
and hence we obtain the following theorem.

\begin{theorem}\label{th61}
If, for some integer $p\geq 1$,  the Nevanlinna function $n$
satisfies one of the assumptions $(j_p)$, $j\in\{1,2,3\}$, and
$D_p\ne 0$ then the following relation holds$:$
\begin{equation}\label{hhf}
n(z)=\left((A-z)^{-1}u,u\right)=-\frac{d_{p-1}(z) n_p(z)+
b_{p-1}d_p(z)} {e_{p-1}(z) n_p(z)+ b_{p-1}e_p(z) },
\end{equation}
where $n_p(z)=\left((D-z)^{-1} u_p, u_p\right)$, $u_p=b_{p-1}e_p$.
\end{theorem}

\begin{remark}\label{npnp-1} (i) Using \eqref{ic} and \eqref{a0b0} we obtain from
\eqref{hhf} with $p=1$:
$$n(z)=-\dfrac{s_0}{z-\dfrac{s_1}{s_0}+n_1(z)}$$
and hence, because $n_1(z)={\rm o}(1)$, $n_1$ is the Schur
transform of $n$: $n_1=\widehat n$. For $p\geq 2$ we obtain from
\eqref{hhf} and \eqref{hhf} with $p$ replaced by $p-1$ and with
the help of \eqref{pol1} and \eqref{lundy} that
$$n_{p-1}(z)=-\dfrac{b_{p-2}^2}{z-a_{p-1}+n_p(z)},$$
hence
$$n_{p-1}(z)=-\dfrac{b_{p-2}^2}{z}-\dfrac{a_{p-1}b_{p-2}^2}{z^2}+{\rm
o}\left(\dfrac{1}{z^2}\right)
$$
and $n_p=\widehat n_{p-1}$.\\
(ii) From \eqref{hhf}, \eqref{rel1}, \eqref{rel2}, and \eqref{AA},
we obtain
$$\widehat n(z)=-\dfrac{s_0}{n(z)}-(z-a_0)=-\frac{\widehat d_{p-2}(z)
n_p(z)+ \widehat b_{p-2}\widehat d_{p-1}(z)} {\widehat e_{p-2}(z)
n_p(z)+ \widehat b_{p-2}\widehat e_{p-1}(z)}.
$$
Since, according to Theorem \ref{orthogpolnhat}, $\widehat d_k$
and $\widehat e_k$ are the polynomials of first and second kind
associated with $\widehat n$ , this formula implies that the
function $n_p$ is the $p-1$-st Schur transform of $\widehat n$.
\end{remark}

The $2\times 2$ matrix polynomial $V$, which generates the
fractional linear transformation \eqref{hhf}, has the property
$$
\det V(z)=b_{p-1}\left( d_p(z)e_{p-1}(z)-e_p(z)d_{p-1}(z)\right)=1.
$$
With
$$
J=\begin{pmatrix}0 & 1
\\-1 & 0\end{pmatrix},
$$
$V$ is $J$-{\it unitary} on the real line, that is,
$$
V(z)JV(z)^*=J \quad z \in \mathbb R.
$$
 Therefore
$V(z)^{-1}$ exists for all $z\in\mathbb C$ and we can form the
polynomial matrix function
$$
V_0(z)=V(z)V(0)^{-1}=
\begin{pmatrix} p_1^{(p)}(z)&p_0^{(p)}(z)\\q_1^{(p)}(z)&q_0^{(p)}(z)
\end{pmatrix}
$$
%\begin{eqnarray*}
%V_0(z)=V(z)V(0)^{-1}&=&\left(\begin{array}{cc}
%-d_{n-1}(z)&-d_n(z)\\e_{n-1} (z)&e_n(z)
%\end{array}
%\right)\left(\begin{array}{cc} -d_{n-1}(0)&-d_n(0)\\e_{n-1}
%(0)&e_n(0)
%\end{array}
%\right)^{-1}\\[2mm]
%&=&\left(\begin{array}{cc}
%p_1^{(n)}(z)&p_0^{(n)}(z)\\q_1^{(n)}(z)&q_0^{(n)}(z)
%\end{array}
%\right)
%\end{eqnarray*}
with
\begin{eqnarray*}
p_0^{(p)}(z)&=&b_{p-1}\left(d_p(z)\,d_{p_-1}(0)-d_{p-1}(z)\,d_p(0)\right),\\
p_1^{(p)}(z)&=&b_{p-1}\left(d_p(z)\,e_{p-1}(0)-d_{p-1}(z)\,e_p(0)\right),
\\
q_0^{(p)}(z)&=&
b_{p-1}\left(e_{p-1}(z)\,d_p(0)-e_p(z)\,d_{p-1}(0)\right),
\\
q_1^{(p)}(z)&=&
b_{p-1}\left(e_{p-1}(z)\,e_p(0)-e_p(z)\,e_{p-1}(0)\right).
\end{eqnarray*}
 A straightforward calculation leads to the relation
$$
n(z)\equiv\left((A-z)^{-1}u,u\right)=\frac{p_1^{(p)}(z)
h_p(z)+p_0^{(p)}(z)} {q_1^{(p)}(z) h_p(z)+q_0^{(p)}(z)},
$$
where
\begin{equation}\label{fraclin}
h_p(z)=-\frac{d_{p-1}(0)n_p(z)+b_{p-1}d_p(0)} {e_{p-1}(0)
n_p(z)+b_{p-1}e_p(0)}.
\end{equation}
With relation \eqref{sc}, the following formulas can be obtained,
compare \cite[I.2.4]{akh}:
\begin{eqnarray*}
p_0^{(p)}(z)&=&z\sum_{k=0}^{p-1}d_k(z)\,d_k(0),\\
p_1^{(p)}(z)&=& 1+ z\sum_{k=0}^{p-1}d_k(z)\,e_k(0),
\\
q_0^{(p)}(z)&=& 1- z\sum_{k=0}^{p-1}e_{k}(z)\,d_{k}(0),
\\
q_1^{(p)}(z)&=&-z\sum_{k=0}^{p-1}e_{k}(z)\,e_{k}(0).
\end{eqnarray*}

\section{Transformation by means of a $u$--resolvent matrix}
\label{uresolvent}

Given again a function $n\in \mathbf N_0$ with one of the
properties $(j_p),\,j=1,2,3$. Besides the decomposition \eqref{uu}
we consider the decomposition
%\begin{equation}\label{vv}
$$
\mathcal H(n)=\mathcal H_{p+1}\oplus\mathcal H'',\qquad \mathcal
H_{p+1}={\rm span}\{\mathcal H_p,\,e_p\}={\rm
span}\{e_0,e_1,\dots,e_p\},
$$
%\end{equation}
and in the space $\calH_{p+1}$ the restriction
$$S:=A|_{\calH_p}=\begin{pmatrix} A_0 \\ B \end{pmatrix}.$$
This restriction is a non-densely defined symmetric operator in
$\calH_{p+1}$ with defect index $(1,1)$, and, evidently, the given
function $n=\left((A-z)^{-1}u,u\right)$ is one of the
$u$--resolvents of this operator $S$. Hence $n$ can be represented
as a fractional linear transformation of some function $g \in
\mathbf N_0$ by means of the $u$--resolvent matrix
$W=(w_{k\ell})_{k,\ell=1}^2$ of $S$:
\begin{equation}\label{c}
n(z)=\frac{w_{11}(z)g(z)+w_{12}(z)}{w_{21}(z)g(z)+w_{22}(z)}.
\end{equation}

Such  a $u$--resolvent matrix $W$ can easily be calculated. To
this end we fix a self-adjoint extension of $S$ in $\calH_{p+1}$,
which means that we fix some $\gamma\in\mathbb R$ in the right
lower corner of the matrix representation of $S$ with respect to
the basis $e_0,e_1,\dots,e_p$ of $\calH_{p+1}$. Denote this matrix
or self-adjoint extension of $S$ in $\calH_{p+1}$ by
$A_{0,\gamma}$:
$$A_{0,\gamma}=\begin{pmatrix} A_0 & \widetilde B  \\ B & \gamma
\end{pmatrix}.$$

 According to \cite{kl2} this $u$--resolvent matrix
$W$ is given by the formula
\begin{equation}\label{gro}
W(z)=\dfrac{1}{(u,\varphi(z^*))}\left( \begin{array}{cc}
(R^\gamma_zu,u)&(R^\gamma_zu,u)Q(z)-(u,\varphi(z^*))(\varphi(z),u)\\
1&Q(z)
\end{array}
\right),
\end{equation}
where $R^\gamma_z=(A_{0,\gamma}-z)^{-1}$, $\varphi(z)$ is a defect
function of $S$ corresponding to the self-adjoint extension
$A_{0,\gamma}$, and $Q$ is the corresponding $Q$--function. An
easy calculation yields
$$
(R^\gamma_zu,u)=\left((A_{0,\gamma}-z)^{-1}u,u\right)=a_{00}(z)-
\frac{b_{p-1}^2a(z)^2}{\Delta(z)},$$ where
$$
\Delta(z)=z-\gamma+b_{p-1}^2a_{11}(z).
$$
Since $S=A_{0,\gamma}|_{\mathcal H_{p}}$ and hence, in terms of
linear relations,
$$S^*=\{\{x,
A_{0,\gamma}x+\lambda e_p\}\,|\, x \in \mathcal H_{p+1}, \lambda
\in \mathbb C\},$$ it is easy to check that for $\varphi(z)$ with
$\{\varphi(z), z\varphi(z)\} \in S^*$ we can choose
$$
 \varphi(z)=(A_{0,\gamma}-z)^{-1} e_{p}
%= \frac 1{\Delta(z)}
%\left(-b_{p-1}(A_0-z)^{-1}e_{n-1}+e_p\right)
= \frac {-1}{\Delta(z)}
 \left(\begin{array}{c}
 -b_{p-1}(A_0-z)^{-1}e_{p-1}\\1
 \end{array}
 \right),
$$
and then the $Q$--function, which is the solution (up to a real
additive constant) of the equation
$$\dfrac{Q(z)-Q(\zeta)^*}{z-\zeta^*}=(\varphi(z),
\varphi(\zeta)),$$ becomes
$$
Q(z)=\frac {-1}{\Delta(z)}.
$$
Inserting these expressions into $W$ from \eqref{gro} we find
\begin{equation}\label{av}
W(z)=\dfrac{\Delta(z)}{b_{p-1}a(z)}\begin{pmatrix}
a_{00}(z)-b_{p-1}^2\dfrac{a(z)^2}{\Delta(z)}&
-\dfrac{a_{00}(z)}{\Delta(z)}\\[4mm]1&\dfrac{-1}{\Delta(z)}
\end{pmatrix}.
\end{equation}
Observe that $W(z)$ is $J$-unitary on the real line. Next we
establish the connection between the matrix functions $V$ from
\eqref{v} and $W$ from \eqref{gro}, in fact we find a simple
expression for $V^{-1}W$. We have
$$
V(z)^{-1}=\dfrac{1}{a(z)} \begin{pmatrix}
-1&a_{00}(z)\\[2mm]
-a_{11}(z)&a_{00}(z)a_{11}(z)-a(z)^2
\end{pmatrix}.
$$
Multiplying this matrix from the right by $W(z)$ from \eqref{av}
 we obtain
 \begin{equation}\label{162}
 V(z)^{-1}W(z)=
- \begin{pmatrix}
 -b_{p-1}&0\\\dfrac{z-\gamma}{b_{p-1}}&-\dfrac{1}{b_{p-1}}
 \end{pmatrix}.
\end{equation}
\begin{theorem}
If the function $n\in \mathbf N_0$ has one of the properties
$(j_p),\,j\in\{1,2,3\}$, then the matrix functions $V$ from
\eqref{v} and $W$ from \eqref{gro} are connected by the relation
\eqref{162}. Therefore for the Nevanlinna functions $n_p$ in
\eqref{le} and $g$ in \eqref{c} the following relation holds$:$
\begin{equation}\label{fg}
n_p(z)=-\frac{b_{p-1}^2}{z-\gamma-\frac 1{g(z)}}.
\end{equation}
\end{theorem}

If $1/g(z)={\rm o}(1)$, then formula \eqref{fg} implies that $n_p$
admits the asymptotic expansion
$$n_p(z)=-\dfrac{b_{p-1}^2}{z}-\dfrac{\gamma b_{p-1}^2}{z^2}+{\rm
o}\left(\dfrac{1}{z^2}\right)$$ and $-1/g$ is its Schur transform:
$-1/g=\widehat n_p$. Hence the number $\gamma$, which defines the
self-adjoint extension of $S$, corresponds to the number $a_{p}$.
\begin{remark}
If instead of a self-adjoint {\it operator} extension $A_{0,\gamma}$ of $S$ we
 choose the (multi-valued) self-adjoint {\it relation} extension of
S:
$$
A_{0,\infty}=S+{\rm span}\,\{0,e_{p}\}=A_0+{\rm
span}\,\{0,e_{p}\},
$$
then we obtain
$$R^\infty_z=(A_{0,\infty}-z)^{-1}=(A_0-z)^{-1}P,$$
where $P$ is the orthogonal projection in $\mathcal H_{p+1}$ onto
$\mathcal H_p$,
$$
\varphi(z)=
\begin{pmatrix} -b_{p-1}(A_0-z)^{-1}e_{p-1} \\
1 \end{pmatrix}, \quad  Q(z)=z+b_{p-1}^2a_{11}(z),
$$
so that
$$
W(z)=-\dfrac{1}{b_{p-1}^2a(z)}
\begin{pmatrix} a_{00}(z) & a_{00}(z+b_{p-1}^2a_{11}(z))-b_{p-1}^2a(z)^2 \\[2mm]
1 & z+b_{p-1}^2a_{11}(z)
\end{pmatrix},
$$
$$
V(z)^{-1}W(z)=\begin{pmatrix}0 & -b_{p-1} \\[2mm]
 \dfrac{1}{b_{p-1}} &
\dfrac{z}{b_{p-1}}
\end{pmatrix},$$ and instead of \eqref{fg} we have
$$
n_p(z)=-\dfrac{b_{p-1}^2}{z+g(z)}.
$$
Thus if $g(z)={\rm o}(1)$, then $n_p$ has the asymptotic expansion
$$n_p(z)=-\dfrac{b_{p-1}^2}{z}+{\rm
o}\left(\dfrac{1}{z^2}\right)$$ and $g$ is the Schur transform of
$n_p$: $g=\widehat n_p$.
\end{remark}

 A more explicit form of the resolvent matrix
$W$ from \eqref{gro} can be obtained following \cite{kl2} and
\cite{abdl}. To this end we decompose the space $\calH_{p+1}$ as
$$
\calH_{p+1}={\rm ran} (S-z)\dotplus {\rm span}\, u,\quad z\in
 \mathbb C,\ a(z)\neq 0,
$$
($\dotplus$ stands for direct sum) and denote for
$y\in\calH_{p+1}$ by ${\textsf{P}}(z)y$ the coefficient of $u$ in
the corresponding decomposition of $y$:
\begin{equation}\label{pots} y=(S-z)x+(\textsf{P}(z)y)\,u
\end{equation}
 with  some $x\in {\rm dom}\,S=\mathcal H_p$.
 Further, define $\textsf{Q}(z)\,y=\left((S-z)^{-1}(y-(\textsf{P}(z)y)u),
 u\right)$. Then, according to \cite{kl2}, the resolvent matrix can be chosen to be
\begin{equation}\label{ww}
W^0(z)=I_2+z\left(\begin{array}{c}\textsf{Q}(z)\\-\textsf{P}(z)\end{array}\right)\left(\textsf{Q}(0)^*\
-\textsf{P}(0)^*\right)J, \quad J=\begin{pmatrix} 0 & 1 \\ -1 & 0
\end{pmatrix}.
\end{equation}

We derive an explicit expression for $W^0(z):=W(z)W(0)^{-1}$,
following \cite{abdl}. To this end, for the vectors and operators
we use matrix representations with respect to the basis $\mathbf
t_0 (\,=u)$, $\mathbf t_1, \ldots, \mathbf t_{p}$.  Recall that
$S_p$ is the Gram matrix associated with this basis. We denote by
$\mathfrak{S}$ the $(p+1)\! \times\! (p+1)$--matrix
$$\mathfrak{S}= \begin{pmatrix} 0&\cdots&\cdots&\cdots&0 \\
     1&0&\cdots&\cdots&0 \\ 0&1&0&\cdots&\vdots \\ \vdots&&\ddots&\ddots&\vdots \\
     0&\cdots&0&1&0
     \end{pmatrix} $$
and by $e$ the first column in the $(p+1)\! \times \!(p+1)$
identity matrix.  Then $S$ and $u$ correspond to $\mathfrak{S}\mid
_{ \mathbb C^p \dotplus \{0\}}$ and $e$. First we apply the
operator $(I-z\mathfrak{S}^*)^{-1}$ to \eqref{pots} and observe
the relation
$$
e^*(I-z\mathfrak{S}^*)^{-1}(\mathfrak{S}-z)x=0,\quad x\in \mathbb
C^p \dotplus \{0\}.
$$
It follows that
\begin{equation}\label{ppp}
e^*(I-z\mathfrak{S}^*)^{-1}y=(\textsf{P}(z)y)\,e^*(I-z\mathfrak{S}^*)^{-1}u
=\textsf{P}(z)y.
\end{equation}
Further, observing that $e^*S_p=
\begin{pmatrix}s_0&s_1&\cdots&s_p\end{pmatrix}$
we obtain
\begin{eqnarray*}
\textsf{Q}(z)y&=&\!\!\!\left((S-z)^{-1}(y-(\textsf{P}(z)y)u),u\right)\\
&=&\!\!\!e^*S_p
\left((\mathfrak{S}-z)^{-1}y-(\mathfrak{S}-z)^{-1}e\left(e^*(I-z\mathfrak{S}^*)^{-1}y\right)
\right)\\
&=&\!\!\!\begin{pmatrix}s_0&s_1&\cdots&s_p\end{pmatrix}
\left((\mathfrak{S}-z)^{-1}(I-z\mathfrak{S}^*)-(\mathfrak{S}-z)^{-1}e
e^*\right)(I-z\mathfrak{S}^*)^{-1}y\\
&=&\!\!\!\begin{pmatrix}s_0&s_1&\cdots&s_p\end{pmatrix}\mathfrak{S}^*(I-z\mathfrak{S}^*)^{-1}y\\
&=&\!\!\!\begin{pmatrix}0&s_0&s_1&\cdots&s_{p-1}\end{pmatrix}(I-z\mathfrak{S}^*)^{-1}y,
\end{eqnarray*}
where for the second last equality sign we have used that
$$
(\mathfrak{S}-z)^{-1}(I-z\mathfrak{S}^*)-(\mathfrak{S}-z)^{-1}ee^*=\mathfrak{S}^*.
$$
Together with \eqref{ppp} we find
$$
\begin{pmatrix}\textsf{Q}(z)\\-\textsf{P}(z)\end{pmatrix}=
\begin{pmatrix}0&s_0&s_1&\cdots&s_{p-1}\\-1&0&0&\cdots&0
\end{pmatrix}\left(I-z\mathfrak{S}^*\right)^{-1},
$$
and \eqref{ww} becomes
$$
W^0(z)=I_2+z\left(\hspace*{-2mm}\begin{array}{cccc}0&s_0&
\cdots&s_{p-1}\\-1&0&\cdots&0\end{array}\hspace*{-2mm}\right)\left(I-z\mathfrak{S}^*
\right)^{-1}S_p^{-1}\left(\hspace*{-2mm}\begin{array}{cc}0
\hspace*{-3mm}&-1
\\s_0\hspace*{-3mm}&0\\
\vdots\hspace*{-3mm}&\vdots\\s_{p-1}\hspace*{-3mm}&0\end{array}
\hspace*{-2mm}\right)\left(\begin{array}{cc}0&1\\-1&0\end{array}\right).
$$

\section{Reproducing kernel spaces: reduction via resolvent invariant subspaces}\label{DQO}

In this section we start from the operator representation of the
Nevanlinna function $n$ in the corresponding reproducing kernel
space $\mathcal L(n)$ with kernel
$$
L_n(z,w)=\dfrac{n(z)-n(w)^*}{z-w^*},\quad z,\zeta\in\mathbb
C\setminus\mathbb R,
$$
see Section \ref{Schurtrsf}, (3). The operator $A$ is introduced
via its resolvent $(A-z)^{-1}$ which is the difference-quotient
operator $R_z$ defined by
\begin{equation}\label{resinv}
(R_zf)(\zeta)=\dfrac{f(\zeta)-f(z)}{\zeta-z}, \quad f \in \mathcal
L(n).
\end{equation}
If $n$ satisfies one of the assumptions $(j_p)$, then, by
\cite[Lemma 5.1]{adlline}, the functions
$$
f_0(\zeta)= n(\zeta),\ f_1(\zeta)= \zeta n(\zeta)+s_0,\ \ldots,\
f_{p}(\zeta)=\zeta^{p} n(\zeta)+\zeta^{p-1}s_0+\cdots +s_{p-1}$$
all belong to ${\mathcal L}(n)$ and
\begin{equation} \label{inp}
\langle f_k, f_j\rangle_{{\mathcal L}(n)}=s_{j+k}, \quad j,\, k =0,1,
\ldots,p.
\end{equation}
In particular, $u:=n\in\mathcal L(n)$, and by the reproducing
property of the kernel $L_n$ we have
$$
n(z)=((A-z)^{-1}u,u)_{{\mathcal L}(n)}.
$$

By $\mathcal U_J$ we denote the class of all $2\times 2$ matrix
polynomials $\Theta$ which are $J$-unitary on $\mathbb R$ and for
which the kernel
$$
K_\Theta(z,w)=\dfrac{J-\Theta(z)J\Theta(w)^*}{z-w^*}
$$
is non-negative. The reproducing kernel Hilbert space with this
kernel will be denoted by $\mathcal H(\Theta)$; its elements are
$2$-vector functions.  The matrix polynomials $V$ and $W$
considered in the previous sections belong to $\mathcal U_J$: this
follows from the Christoffel--Darboux formulas \eqref{sc} for $V$
and from \eqref{162} for $W$. Note that if $\Theta$ belongs to
$\mathcal U_J$, then ${\rm det}\, \Theta(z)\equiv c$, where is $c$
is a unimodular complex number, because the determinant ${\rm
det}\,\Theta(z)$ is a non-vanishing polynomial in $z$.

The following theorem was proved in \cite[Theorem 8.1]{adls}, even
in an indefinite setting.

\begin{theorem} \label{rkps1}
Let\ $n\in\mathbf N_0$ and suppose that there exists a matrix
polynomial
$$
\Theta=\begin{pmatrix} a & b \\ c & d \end{pmatrix} \in\mathcal
U_\ell
 $$
 such that the mapping
$$
\mathbf u \longrightarrow \begin{pmatrix}1 &-n
\end{pmatrix}\mathbf u
$$
is an isometry from $\mathcal H(\Theta)$ into $\mathcal L(n)$.
Define the function $\check n$ by
$$
n(z)=\dfrac{a(z)\check n(z)+b(z)}{c(z)\check n(z)+d(z)}.
$$
Then the following statements hold.
\begin{enumerate}
\item[\rm (i)] $\check n$ is Nevanlinna function.\vspace*{0.5mm}

 \item[\rm (ii)] The
mapping $g \mapsto f$:
$$
f(\zeta)=\big(a(\zeta)-n(\zeta)c(\zeta)\big)g(\zeta)
$$
is an isometry from $\mathcal L(\check n)$ into $\mathcal L(n)$.\vspace*{0.5mm}
\item[\rm (iii)] We have
$$
\mathcal L(n)=\begin{pmatrix}1& -n
\end{pmatrix}\mathcal H(\Theta)\oplus(a-nc)\mathcal
L(\check n)
$$
and the mapping
$$
W: \mathcal L(n)\ni
f\mapsto \begin{pmatrix} \mathbf u \\
g \end{pmatrix} \in\begin{pmatrix}\mathcal H(\Theta)\\
\mathcal L(\check n)\end{pmatrix},
$$
where $f, \mathbf u$, and $g$ are connected by the relation
%\begin{equation}\label{Sorbonne}
$$
f(\zeta)=\begin{pmatrix}1 &- n(\zeta)\end{pmatrix}\mathbf
u(\zeta)+\big(a(\zeta)-n(\zeta)c(\zeta)\big)g(\zeta),
$$
%\end{equation}
is a unitary mapping from $\mathcal L(n)$ onto $\mathcal
H(\Theta)\oplus\mathcal L(\check n)$.\vspace*{0.5mm}

\item[\rm (iv)] The mapping $WR_zW^*$ is of the form
\begin{equation}\label{henri4}
 WR_zW^*=\begin{pmatrix} R_{11}(z) & R_{12}(z) \\R_{21}(z) & R_{22}(z)
 \end{pmatrix}: \begin{pmatrix} \mathcal P(\Theta) \\ \mathcal L(\check n)
 \end{pmatrix}
 \rightarrow \begin{pmatrix} \mathcal P(\Theta) \\ \mathcal L(\check n)
 \end{pmatrix},
 \end{equation}
with
$$
\begin{array}{ll}
R_{11}(z)&=R_z-\dfrac{1}{k(z)} (R_z\Theta)(\,\cdot\,)
\begin{pmatrix} \check n(z) \\1 \end{pmatrix}
\begin{pmatrix} 0 & 1 \end{pmatrix}
 \,E_z \\[3mm]
       &=R_z-K_{\Theta}(\,\cdot\,,z^*)\begin{pmatrix}
1 \\ -n(z) \end{pmatrix}\begin{pmatrix} 0 & 1 \end{pmatrix}\, E_z,
\\[3mm]
R_{12}(z)&= \dfrac{1}{k(z)}(R_z\Theta)(\,\cdot\,) \begin{pmatrix} d(z) \\[1mm]
-c(z)
\end{pmatrix} \,E_z \\
 & =-(a(z)-n(z)c(z))K_{\Theta}(\,\cdot\,,z^*)
\begin{pmatrix} 0 \\ 1 \end{pmatrix}\, E_z,
\\[3mm]
R_{21}(z)&=-\dfrac{1}{k(z)}(R_z\check n)(\,\cdot\,)\begin{pmatrix}
0 & 1
\end{pmatrix}\, E_z \\[3mm]
&=-\dfrac{1}{k(z)} L_{\check n}(\,\cdot\,,z^*)\begin{pmatrix} 0 &
1
\end{pmatrix}\, E_z,
\\
R_{22}(z)&= R_z - \dfrac{c(z)}{k(z)}(R_z\check n)(\,\cdot\,)\,E_z \\[3mm]
 &=R_z - \dfrac{c(z)}{k(z)}L_{\check n}(\,\cdot\,,z^*)\,E_z,
\end{array}
$$
where $R_z$  is the difference-quotient operator, $E_z$ is the
operator of evaluation at the point $z$ on any reproducing kernel
space, and

%\begin{equation}\label{k}
$$
k(z)=c(z)\check n(z)+d(z)=\dfrac{{\rm
det}\,\Theta(z)}{a(z)-n(z)c(z)}.
$$
%\end{equation}
\end{enumerate}
\end{theorem}
We mention that formula \eqref{henri4} corresponds to the relation
\eqref{respp} above.

A space of functions is called {\it resolvent-invariant} if it is
invariant under the difference-quotient operator $R_z$ as defined
in \eqref{resinv}. In the following lemma, with a
resolvent-invariant non-degenerate invariant subspace of a certain
inner product space a  $2\times 2$ matrix function is associated.

\begin{lemma}
Let ${\mathcal M}$ be a finite dimensional resolvent--invariant
space of $2$--vector polynomials endowed with an inner product
$\langle\, \cdot\, ,\,\cdot\, \rangle$ such that
\begin{equation}
\label{eq} \langle R_zf,g\rangle-\langle f,R_wg\rangle-
(z-w^*)\langle R_zf,R_wg\rangle=g(w)^*Jf(z),
\end{equation}
and let that ${\mathcal M}_1$ be a resolvent--invariant
non-degenerate subspace of ${\mathcal M}$. Then there exists a
$\Theta_1\in{\mathcal U_J}$ such that
\begin{enumerate}
\item[\rm (i)] ${\mathcal M}_1={\mathcal
H}(\Theta_1)$,\vspace*{0.5mm} \item[\rm (ii)]  ${\mathcal
M}={\mathcal H}(\Theta_1)\oplus \Theta_1{\mathcal N} $ where
${\mathcal N}=\Theta_1^{-1}{\mathcal M}_1^{\perp}$ is a
resolvent--invariant space of $2$-vector polynomials, for which
the relation \eqref{eq} holds if equipped with the inner product
$$
( \Theta_1^{-1}f, \Theta_1^{-1}g)_{{\mathcal N}}= \langle f,
g\rangle,\quad f,g\in{\mathcal M_1}^\perp.
$$ \label{ad}
\end{enumerate}
\end{lemma}

Relation \eqref{eq} is often  called {\it de Branges identity},
see \cite{de-hsafI} and,  for further references, \cite{dym2}.
That ${\mathcal N}$ consists of $2$--vector polynomials is due to
fact that $\Theta_1^{-1}(z)=-J\Theta_1(z^*)^*J$ is a matrix
polynomial. The other claims  of the lemma follow from
\cite[Theorem 3.1]{ad-90}.

Now we formulate and prove Theorem \ref{th61} again in the context
of reproducing kernel spaces.

\begin{theorem}\label{metro}
If, for some integer $p\geq 1$,  the Nevanlinna function $n$
satisfies one of the assumptions $(j_p)$, $j\in\{1,2,3\}$, and
$D_p\ne 0$ then the following relation holds$:$
\begin{equation}\label{hhfagain}
n(z)=-\frac{d_{p-1}(z) \check n_p(z)+ b_{p-1}d_p(z)} {e_{p-1}(z)
\check n_p(z)+ b_{p-1}e_p(z) }.
\end{equation}
where  $\check n_p$ is a Nevanlinna function such that $\check
n_p\in {\mathcal L}(\check n_p)$, and
\begin{equation}\label{88}
\check n_p(z)=(R_z\check n_p,\check n_p)_{{\mathcal L}(\check
n_p)}.
\end{equation}
\end{theorem}

Comparing \eqref{hhf} and \eqref{hhfagain} we find that $\check
n_p(z)=n_p(z)$, the $p$-th element in the sequence obtained by
applying the Schur transformation $p$ times starting with $n$.

\begin{proof}[Proof of Theorem {\rm \ref{metro}}]
Let $\mathcal M$ be the linear space spanned by the $p+1$
$2$-vector functions
\begin{equation}\label{mbasis}
\mathbf f_0(\zeta)= \begin{pmatrix} 0 \\ -1 \end{pmatrix}, \mathbf
f_1(\zeta)=\begin{pmatrix} s_0 \\ -\zeta \end{pmatrix},
\dots, \mathbf f_p(\zeta)=\begin{pmatrix} s_0 \zeta^{p-1}+ \cdots+ s_{p-1} \\
-\zeta^{p}
\end{pmatrix}
\end{equation}
and equipped with the inner product which makes the map $\mathbf u
\mapsto
\begin{pmatrix}
 1 & -n
 \end{pmatrix} \mathbf u $ an isometry from ${\mathcal
M}$ into ${\mathcal L}(n)$, see \eqref{inp}. Note that if $n$ has
the integral representation \eqref{int}, the elements of $\mathcal
M$ are of the form
$$\begin{pmatrix} \int_{-\infty}^\infty (R_\zeta
f)(t)\,d\sigma(t) \\ -f(\zeta) \end{pmatrix},
$$
where $f$ is a polynomial of degree $\leq p$. Indeed, it suffices
to show this for the basis elements of $\mathcal M$: If
$f(\zeta)=\zeta^j$, then
$$R_\zeta f(t)=\dfrac{\zeta^j-t^j}{\zeta-t}=\zeta^{j-1}+t\zeta^{j-2}+\cdots
+t^{j-2}\zeta+t^{j-1},$$ and hence, on account of \eqref{Achmom},
$$
\int_{-\infty}^\infty (R_\zeta f)(t)\,d\sigma(t)= s_0
\zeta^{j-1}+s_1\zeta^{j-2}+\cdots+s_{p-2}\zeta+s_{p-1}.
$$
It follows that $\mathcal M$ is also spanned by the polynomial
vectors
$$\begin{pmatrix} -d_j \\ e_j \end{pmatrix}, \quad
j=0,1,\ldots,p,$$ where $e_j$ and $d_j$ are the polynomials of
first and second kind associated with $n$, see \eqref{innere} and
\eqref{innerd}.

Let ${\mathcal M}_p$ be the space spanned by the first $p$ of the
$2$-vector functions in \eqref{mbasis}. Since $S_{p-1}$ from
\eqref{hankel} is a positive matrix, the space ${\mathcal M}_p$ is
non-degenerate and
%\begin{equation}\label{AlpDym}
$$
\mathcal M=\mathcal M_p \oplus {\rm span} \begin{pmatrix}- d_p \\
e_p \end{pmatrix}.
$$
%\end{equation}
 As both $\mathcal M$ and
$\mathcal M_p$ are resolvent-invariant spaces, by Lemma \ref{ad}
we have that for some $\Theta_1\in{\mathcal U_J}$, which is
normalized by $\Theta_1(0)=I_2$ (and hence ${\rm det}\,\Theta(z)
\equiv 1$),
%\begin{equation}\label{ADdecomp}
$$
{\mathcal M}_p={\mathcal H}(\Theta_1), \quad {\mathcal
M}={\mathcal H}(\Theta_1)\oplus \Theta_1{\mathcal N}.
$$
%\end{equation}
Here ${\mathcal N}$ is a one-dimensional resolvent-invariant
space, which, when equipped with the induced inner product,
satisfies the de Branges identity and therefore is spanned by a
constant $J$--neutral vector $\begin{pmatrix} \alpha & \beta
\end{pmatrix}^\top$ such that
%\begin{equation}\label{ab}
$$
\Theta_1(z) \begin{pmatrix} \alpha \\ \beta \end{pmatrix}=
b_{p-1}\begin{pmatrix} -d_p(z) \\ e_p(z) \end{pmatrix}.
$$
%\end{equation}
For $\lambda \in \mathbb R$ denote by $C_{\lambda}$ the constant
$J$--unitary matrix
$$C_\lambda=\left\{ \begin{array}{ll}
\begin{pmatrix} \lambda \alpha & \alpha \\ -\alpha^{-*}
+\lambda \beta & \beta \end{pmatrix}, & \alpha \neq 0,\\[4mm]
\begin{pmatrix} \beta^{-*} & 0 \\ \lambda \beta & \beta \end{pmatrix},
& \alpha = 0.\\
\end{array}\right.$$
Then there exists a $\lambda$ such that
$$\Theta(z):=\Theta_1(z)C_\lambda =\begin{pmatrix}
a(z)&-b_{p-1}d_p(z)\\ c(z)& b_{p-1} e_p(z)\end{pmatrix},
$$
where $a$ and $c$ are polynomials such that $\deg c < \deg e_p=p$.
The inclusion
$$ R_0 \begin{pmatrix} a \\ c \end{pmatrix} \in {\mathcal
H}(\Theta) ={\mathcal H}(\Theta_1)=\mathcal M_p
$$
implies that $\deg a <p-1$.  From $\det \Theta(z) \equiv 1$ it
follows that $x=a$ and $y=c$ are polynomial solutions of the
equation
$$
x(z)e_p(z)+y(z)d_p(z)= \dfrac{1}{b_{p-1}}.
$$
Since all polynomial solutions of this equation are given by
$$
x(z)=a(z)-s(z)d_p(z), \quad y(z)=c(z)+s(z)e_p(z)
$$
with some polynomial $s$, the solutions $x=a$ and $y=c$ have
minimal degrees and because of that they are unique. Observing
\eqref{lundy}, we find
$$a(z)=-d_{p-1}(z), \quad c(z)=e_{p-1}(z).$$
Hence
$$
\Theta(z)=\begin{pmatrix} -d_{p-1}(z) &-b_{p-1}d_p(z)\\
e_{p-1}(z)& b_{p-1}e_p(z)\end{pmatrix}=V(z)
$$
and $C_\lambda=\Theta(0)$ is the coefficient matrix of the
fractional linear transformation \eqref{fraclin}.

Define the function $\check n_p$ by \eqref{hhfagain}. Then,
according to Theorem \ref{rkps1}, it is a Nevanlinna function. We
show that $ \check n_p \in \mathcal L(\check n_p)$. The function
$f_p(\zeta)=\begin{pmatrix} 1 & -n(\zeta) \end{pmatrix} \mathbf
f_p(\zeta)$ belongs to $\mathcal L(n)$ and, according to Theorem
\ref{rkps1} (iii), it can be written as
$$
\begin{pmatrix} 1 & -n(\zeta) \end{pmatrix} \mathbf f_p(\zeta)
=\begin{pmatrix} 1 & -n(\zeta) \end{pmatrix} \mathbf u_p(\zeta)+
(a(\zeta)-n(\zeta)c(\zeta))g_p(\zeta)
$$
with $\mathbf u_p \in \mathcal M_p$, $g_p \in \mathcal L(\check
n_p)$, and the two summands on the righthand side are orthogonal.
This orthogonality and the isometry of the mapping
$\begin{pmatrix} 1 & -n \end{pmatrix}$ imply that $(0 \neq)\,
\mathbf f_p-\mathbf u_p \in \mathcal M_p^\perp$ and hence there is
a non-zero complex number $\gamma$ such that
$$\mathbf f_p-\mathbf u_p  =\gamma \begin{pmatrix} -d_p(\zeta)  \\ e_p(\zeta)
\end{pmatrix}.
$$
Therefore
$$(a(\zeta)-n(\zeta)c(\zeta))g_p(\zeta)=
-\gamma \left(d_p(\zeta)+e_p(\zeta)n(\zeta)\right)
$$
and
$$g_p(\zeta)=-\gamma
\dfrac{e_p(\zeta)\,n(\zeta)+d_p(\zeta)}{-n(\zeta)\,c(\zeta)+a(\zeta)}
=-\dfrac{\gamma}{b_{p-1}} \check n_p(\zeta).$$ Hence $\check n_p
\in \mathcal L(\check n_p)$. Equality \eqref{88} follows from item
(3) in Section \ref{saopreal}.
\end{proof}

\end{document}